\documentclass{article}
\usepackage[a4paper]{geometry}
\usepackage{amssymb,amsmath,amsthm,verbatim,IEEEtrantools,macros}
\author{Robert Fraser}
\title{A framework for constructing sets without configurations}
\begin{document}
\maketitle
\begin{abstract}
We discuss a framework for constructing large subsets of $\mathbb{R}^n$ and $K^n$ for non-archimedean local fields $K$. This framework is applied to obtain new estimates for the Hausdorff dimension of angle-avoiding sets and to provide a counterexample to a limiting version of the Capset problem.
\end{abstract}
\section{Introduction and Background}
Many questions in additive combinatorics and geometric measure theory are of the following form: If a set $S$ in some space $X$ is large in an appropriate sense, then must it contain a certain configuration of points? The techniques involved in studying the problem depend upon the space $X$ and the configuration of points being studied. Problems in additive combinatorics are often concerned with the case in which $X$ is a finite abelian group, $S$ is assumed to contain a certain number of elements depending on the order of $X$, and the configurations being studied are solutions to linear equations in $X$. For example, Roth's theorem on $3$-term arithmetic progressions \cite{Roth52} and the recent capset result of Ellenberg and Gijswijt \cite{EllenbergGijswijt17} are of this type. In geometric measure theory, the space $X$ is often taken to be $\mathbb{R}^n$ and the configurations under study tend to be geometric in nature, and the results concern the Hausdorff dimension of the set $S$. Examples include the recent result of Harangi et al \cite{HarangiKeletiKissMagaMatheMattilaStrenner13} on angle-avoiding sets and the general work of Andr\'as M\'ath\'e on polynomial configurations \cite{Mathe17}.

Given a commutative ring $R$ and a function $f : R^{nv} \to R$, we are interested in subsets of $R^n$ with large Hausdorff dimension not containing any $v$ distinct points $x_1, \ldots, x_v$ such that $f(x_1, x_2, \ldots, x_v) = 0$. M\'ath\'e \cite{Mathe17} considers the case in which $R = \mathbb{R}$ and $f$ is a polynomial of degree $d$ with rational coefficients, obtaining a Hausdorff dimension bound of $n/d$. In particular, this bound does not depend on the number of points $v$ in the configuration. M\'ath\'e applies the $n/d$ bound to obtain a result on angle-avoiding sets. The author and Pramanik \cite{FraserPramanik18} obtain a bound of $\frac{1}{v-1}$ for non-polynomial functions $f$ satisfying some mild conditions on the derivatives.

We will obtain bounds for functions $f$ admitting a special set of points called a \textbf{landmark pair}. A landmark pair is a ubiquitous set of points that avoid a neighbourhood of $0$ and satisfy certain mapping properties under $f$. The main result of this paper is Theorem \ref{LandmarkTheorem}, which allows for the construction of sets $E$ of large Hausdorff dimension avoiding such functions $f$. This theorem implies a slight generalization of M\'ath\'e's result:

\begin{mycor}\label{MatheAlgebraicResult}
Let $p_j(x_1, \ldots, x_{v_j}) : \mathbb{R}^{nv_j} \to \mathbb{R}$ be a countable collection of polynomials of degree at most $d$  whose coefficients are algebraic over the rational numbers. Then there exists a subset $E \subset \mathbb{R}^n$ with Hausdorff dimension $\frac{n}{d}$ that does not contain, for any $j$, any $v_j$-tuple of distinct points $x_1, \ldots, x_{v_j}$ such that $p_j(x_1, \ldots, x_{v_j}) = 0$.
\end{mycor}

Theorem \ref{LandmarkTheorem} can be applied to a diverse set of avoidance problems on a variety of spaces. For example, we are also able to obtain a $p$-adic version of M\'ath\'e's result (Corollary \ref{LocalFieldPolynomial}) using the main theorem in this paper.
\section{Landmark Systems}
The main result is the following:
\begin{mythm}\label{LandmarkTheorem}
Let $K$ be either a non-archimedean local field or $\mathbb{R}$, and let $\{ f_{q} \}_{q = 1}^{\infty} : K^{nv_{q}} \to K$ be a sequence of $|\alpha_q|$-times strictly differentiable functions. Suppose that there exists a ball $B \subset K^n$ on which each function $f_q$ has some partial derivative $\partial^{\alpha_{q}}$ of order $|\alpha_{q}|$ that does not vanish on $B^{v}$ for any $q$. Suppose there exists a sequence $\epsilon_m$ of positive real numbers with limit $0$, and a sequence of weak approximate landmark pairs $\{(\ell_1^{(q, m, j)}, \ell_2^{(q, m, j)})\}_{q \in \Z, m \in \Z, j \in \{0, 1, \ldots, |\alpha_q| - 1\}}$ with parameters $r, \gamma, \sigma$ and of degree $d + \epsilon_m$, and adapted to $D_{j, q} f_q$, where $\{D_{j, q}\}_{j=0}^{|\alpha_q|} $ is a sequence of differential operators such that $D_{0, q}$ is the identity operator, $D_{j, q} = D_{j-1, q} \partial_{x_{i_0(j,q)}^{(k_0(j, q))}}$ and such that $D_{|\alpha_q|, q} = \partial^{\alpha_q}$.  Then there exists a set $E \subset B$ of Hausdorff dimension $\frac{n \sigma}{d \gamma}$ that such that $E$ does not contain any $v_q$ distinct points $(x^{(1)}, \ldots, x^{(v_q)})$ such that $f_q(x^{(1)}, \ldots, x^{(v_q)}) = 0$ for any $q$.
\end{mythm}
In order to make sense of this result, we need to define the notion of a weak approximate landmark pair.

Let $R$ be a commutative ring equipped with a metric $\rho(x,y) : R \times R \to \mathbb{R}$ satisfying the following properties for all $x, y$, and $z$

\begin{IEEEeqnarray}{rCl}
\rho(x + z, y + z) & = & \rho(x,y) \\
\rho(0, x^n) & = & \rho(0,x)^n
\end{IEEEeqnarray}

Suppose that, with respect to this metric, $R$ is locally compact and does not have any isolated points.

\begin{mydef}[Landmark System]\label{LandmarkSystem} We will call $\{\ell_w : R \to \mathbb{Z}_+ \cup \infty\}_{w=1}^{\infty} $ a landmark system for addition and multiplication on $\Omega$; $\Omega \subset R$ compact, if it satisfies the following properties for some positive real numbers $r$, $\gamma$, and $\sigma$ such that $\gamma \geq \sigma$, and some appropriate $\phi_1$ and $\phi_2$:
\begin{itemize}
\item Monotonicity property: $\ell_w(x) \leq \ell_v(x)$ whenever $w > v$ and $x \in \Omega$
\item Additive property: \[\ell_{\phi_1(w_1, w_2)} (x + y) \leq \max(\ell_{w_1}(x), \ell_{w_2}(y)) + o(\ell_{w_1}(x)+ \ell_{w_2}(y)) \] for any $x, y \in \Omega$.
\item Multiplicative property: \[\ell_{\phi_2(w_1,w_2)}(xy) \leq \ell_{w_1}(x) + \ell_{w_2}(y) + o(\ell_{w_1}(x) + \ell_{w_2}(y)) \] for any $x, y \in \Omega$.
\item Separation property: The ball $B(0; C_1(\epsilon, w)^{-1} r^{\gamma j (1 + \epsilon)})$ does not contain any nonzero points $y$ such that $\ell_w(y) \leq j$.  In particular, this will hold if the points satisfying $\ell_w(y) \leq j$ are $C_1(\epsilon, w)^{-1} r^{ \gamma j (1 + \epsilon)}$-separated.
\item Ubiquity property: For any $x \in \Omega$ and any integer $k > 0$, and any $\epsilon > 0$ there is at least one point $y$ in the ball  $B(x;C_2(\epsilon,w) r^{(\sigma  - \epsilon) k})$ such that $\ell_w(y) \leq k$. 
\end{itemize}
\end{mydef}
If there exists a function $\ell(x)$ such that $\ell_w(x) = \ell(x)$ for all $w$ and all $x$, we call $\ell$ a landmark function for $R$. In practice, the error terms in the additive and multiplicative property will be very small; in every example we present, they can in fact be taken to be constant or zero.

We will not be particularly concerned with landmark systems per se, but will instead concern ourselves with the related notion of a landmark pair.
\begin{mydef}[Landmark Pair]\label{LandmarkPairDef}
Let $X$ be a locally compact metric space, let $Y$ be a pointed metric space with distinguished point $0$ and let $f : X^v \to Y$ be a function. For positive real numbers $r, \gamma, \sigma$ with $\gamma \geq \sigma$, we will call $(\ell_1, \ell_2)$, where $\ell_1 : X \to \mathbb{Z}_+ \cup \{\infty\},$ $\ell_2 : Y \to \mathbb{Z}_+ \cup \{\infty\}$ an $(r, \gamma, \sigma)$-landmark pair adapted to $f$ of degree $d$ on a compact set $\Omega \subset X$ if $(\ell_1, \ell_2)$ satisfies the following properties for all $\epsilon > 0$:
\begin{itemize}
\item Function property: 
\[\ell_2(f(x_1, \ldots, x_v)) \leq d \max_{1 \leq i \leq v} \ell_1(x_i) + o(\max_{1 \leq i \leq v} \ell_1(x_i))\]
for all $x_1, \ldots, x_v \in \Omega$.
\item Separation property: There are no nonzero points $y$ in $B(0; C_1(\epsilon)^{-1} r^{\gamma j(1 + \epsilon)})$ such that $ \ell_2(y) \leq j$. In particular, this will hold if the points satisfying $\ell_2(y) \leq j$ are $C_1(\epsilon)^{-1} r^{\gamma j (1 + \epsilon)}$-separated.
\item Ubiquity property: For any $x \in \Omega$, any integer $k > 0$, and any $\epsilon > 0$ there is at least one $y$ in the ball $B(x;C_2(\epsilon)r^{(\sigma - \epsilon) k })$ such that $\ell_1(y) \leq k$.
\end{itemize}
\end{mydef}
In the context of a landmark pair $(\ell_1, \ell_2)$, points for which $\ell_1$ or $\ell_2$ is finite will be called \textbf{landmarks}.
For us, the primary interest in landmark systems is that they give rise to landmark pairs for polynomials with coefficients in the landmark system. 
\begin{mylem}\label{PolynomialLandmarkLemma}
If $\{\ell_w\}_{w=1}^{\infty}$ is an $(r, \gamma, \sigma)$-landmark system on $\Omega \subset R$, and if $p$ is a polynomial of degree $d$ whose coefficients are finite with respect to some $\ell_{w_1}$, then there exists a $w(p)$ such that $(\ell_1, \ell_{w(p)})$ is an $(r, \gamma, \sigma)$-landmark pair on $\Omega$ of degree $d$ for $p$. 
\end{mylem}
\begin{proof}
The separation and ubiquity conditions immediately follow for any $w > 1$ from the assumption. We verify the function condition by induction on the degree of the polynomial. We begin with polynomials of degree $1$. 

In order to show the statement for polynomials of degree $1$, we will first suppose $p$ is of the form $p(x_1, \ldots, x_v) = a_1 x_1 + a$ for some $a, a_1$ such that $\ell_{w_1}(a)$ and $\ell_{w_1}(a_1)$ are finite. Then $\ell_{\phi_1(\phi_2(w_1, 1), w_1)}(a_1 x + a)$ can be estimated by successively applying the additive and multiplicative conditions:
\begin{IEEEeqnarray*}{Cl}
& \ell_{\phi_1(\phi_2(w_1, 1), w_1)}(a_1 x + a) \\
\leq & (1 + o(1)) \max(\ell_{\phi_2(w_1, 1)}(a_1 x), \ell_{w_1}(a))\\
\leq & (1 + o(1)) \max((1 + o(1))(\ell_{w_1}(a_1) + \ell_1(x)) , \ell_{w_1}(a))\\
\leq & (1 + o(1)) \ell_1(x) + \ell_{w_1}(a_1) + \ell_{w_1}(a) \\
\leq & (1 + o(1)) \ell_1(x)
\end{IEEEeqnarray*}
Thus the function condition holds for $p$.

Next, we will show the statement for arbitrary linear polynomials by performing an induction on the number of linear terms. Suppose that we know the condition holds for polynomials of the form $a_1 x_1 + \cdots + a_{v-1} x_{v-1} + a$. We will show that the condition holds for polynomials of the form $a_1 x_1 + \cdots + a_v x_v + a$. Let $p(x_1, \ldots, x_v) = a_1 x_1 + a_2 x_2 + \cdots + a_v x_v + a$, and let $r(x_2, \ldots, x_v) = a_2 x_2 + \cdots + a_v x_v + a$, so that $p(x_1, \ldots, x_v) = a_1 x_1 + r(x_2, \ldots, x_v)$. Then 
\begin{IEEEeqnarray*}{Cl}
& \ell_{\phi_1(\phi_2(w_1,1),w(r))}(a_1 x_1 + r(x_2, \ldots, x_v)) \\
\leq & (1 + o(1)) \max(\ell_{\phi_2(w_1, 1)}(a_1 x_1), \ell_{w(r)}r(x_1, \ldots, x_v)) \\
\leq & (1 + o(1)) \max( \ell_{w_1}(a_1) + \ell_1 (x_1) + o(\ell_1(x_1)), \max_{j \geq 2}(\ell_1(x_j)) + o(\sum_{j \geq 2}l_1(x_j))) \\
\leq & (1 + o(1)) \max(\ell_1(x_1), \ldots, \ell_1(x_v)).
\end{IEEEeqnarray*}
This proves the statement for polynomials of degree $1$.

Now, suppose we know the statement is true for polynomials of degree $d - 1$. We will show it is true for polynomials of degree $d$. 

To show this, we will induct on the number of terms $t$ of degree $d$. If $t = 1$, then $p(\mathbf{x}) = a_{\alpha} \mathbf{x}^{\alpha} + q(\mathbf{x})$, where $\alpha$ is some multi-index of degree $d$ and $q(x)$ has degree at most $d - 1$.  Let $w^{***}$ be the value of $w(r)$ corresponding to the polynomial $r(x) = x^{\beta}$, where $\beta$ is a multi-index of degree $d-1$ that is obtained from $\alpha$ by decrementing the first nonzero entry of $\alpha$; say the $x_i$ entry.  Let $w^{**}$ be $\phi_2(1, w^{***})$, let $w^* = \phi_2(w^{**}, w_1)$, and let $w = \phi_1(w^*, w(q))$. We claim that $w(p)$ can be chosen to be $w$. To see this, let $\ell_1(\mathbf{x}) = \max_j \ell_1(x_j)$ and observe
\begin{IEEEeqnarray*}{Cl}
& \ell_{\phi_1(w^*, w(q))}(a_{\alpha} \mathbf{x}^{\alpha} + q(\mathbf{x})) \\
\leq & (1 + o(1))\max(\ell_{w^*}(a_{\alpha} \mathbf{x}^{\alpha}), \ell_{w(q)}(q(\mathbf{x})))  \\
\leq & (1 + o(1)) \max((1 + o(1)) (\ell_{w_1}(a_{\alpha}) + \ell_{w^{**}}(\mathbf{x}^{\alpha})), (1 + o(1))(n-1)\ell_1(\mathbf{x}))  \\
\leq & (1 + o(1))\max((1 + o(1))(\ell_{w_1}(a_n) + \ell_{1}(x_i) + \ell_{w^{***}}(\mathbf{x}^{\beta})), (n-1) \ell_1(\mathbf{x})) \\
\leq & (1 + o(1))\max(\ell_{w_1}(a_n) + \ell_1(\mathbf{x}) + (n-1) \ell_1(\mathbf{x})), (n-1) \ell_1(\mathbf{x})))
\end{IEEEeqnarray*}
and this maximum is no more than $n \ell_1(\mathbf{x}) + o(\ell_1(\mathbf{x}))$, as desired.

Finally, if $t > 1$, then we can write $p(x) = a_{\alpha} x^{\alpha} + q(x)$, where $q(x)$ is a polynomial of degree $d$ with $t-1$ terms of degree $d$.  A similar argument to the inductive step in the linear case above gives the desired result.
\end{proof}
In fact, for our purposes, we do not need landmarks to map exactly to other landmarks, but only to map to points that are close to landmarks. We codify this notion in the following definition.
\begin{mydef}[Approximate Landmark Pair]\label{ApproximateLandmarkPairDef}
Let $X$ be a locally compact metric space, $\Omega \subset X$ compact, and $Y$ a pointed metric space with distinguished point $0$. Let $r, \gamma$, and $\sigma$ be positive real numbers such that $\gamma \geq \sigma$. A \textbf{weak approximate landmark pair of degree $d$ adapted to $f$ on $\Omega$} is a pair of functions $(\ell_1, \ell_2)$ satisfying the separation and ubiquity conditions for $\gamma$ and $\sigma$ such that there exists an infinite subset $J \subset \mathbb{N}$ and a number $\epsilon > 0$ such that, for any $j \in J$ and any $x_1, \ldots, x_v$ satisfying $\max(\ell_1(x_1), \ldots, \ell_1(x_v)) = j$, we have that
\[\ell_2(y) \leq d j + o(j)\]
for some $y$ satisfying the condition that $\rho(y, f(x)) \leq r^{(\gamma d + \epsilon)j}$.
If $J = \mathbb{N}$, we call $(\ell_1, \ell_2)$ an \textbf{approximate landmark pair of degree $d$ adapted to $f$}.
\end{mydef}
We provide some examples of landmark systems, landmark pairs, and weak approximate landmark pairs in the following sections.
\section{Examples of Real Landmarks}
We begin with a motivating example. This example comes from \cite{Mathe17} and serves as the motivation for landmark systems.
\begin{myex}\label{RationalExample}
Let $R = \mathbb{R}$. Let $N$ be a fixed integer. Define $\ell(x)$ to be the minimal nonnegative integer $n$ (if such an $n$ exists) such that 
\[x = \frac{a}{N^n}\]
where $a$ is an integer. Take $\ell(x) = \ell_w(x)$ for all $w$. Then $\{\ell_w\}_{w=1}^{\infty}$ is a landmark system for $\mathbb{R}$ with $r = N^{-1}$, $\gamma = 1$, and $\sigma = 1$.
\end{myex}
\begin{proof}
The monotonicity property is trivially satisfied because $\ell_w$ does not depend on $w$.

Consider the sum $\frac{a}{N^{n_1}} + \frac{b}{N^{n_2}}$. Without loss of generality, suppose $n_2 \geq n_1$. Then we can rewrite the sum as $\frac{a N^{n_2 - n_1} + b}{N^{n_2}}$, and therefore $\ell_w$ satisfies the additive property (the function $\phi_1$ is not important because $\ell_w$ does not depend on $w$).

The product $\frac{a}{N^{n_1}} \cdot \frac{b}{N^{n_2}}$ is equal to $\frac{ab}{N^{n_1 + n_2}}$, so the multiplicative property is satisfied.

Clearly, the multiples of $\frac{1}{N^n}$ are $\frac{1}{N^n}$-separated, so the separation condition is satisfied for $\gamma = 1$.

Finally, each half-open ball of radius $\frac{1}{N^n}$ contains a number of the form $\frac{a}{N^n}$ for some integer $a$, so the ubiquity condition is satisfied for $\sigma = 1$.
\end{proof}
The above example is somewhat trivial in that the $\epsilon$ and $w$ from the definition were not necessary. We present a more nontrivial example to illustrate the purpose of the $w$ and $\epsilon$ in the definition.
\begin{myex}\label{AlgebraicExample}
Let $R = \mathbb{R}$. Let $\mathbb{Q}(\theta)$ be a finite real extension of $\mathbb{Q}$ of degree $k$ and let $\ell_w(x)$ be the minimal $n$ (if it exists) such that
\[x = \frac{1}{2^n} (a_0 + \cdots + a_{k-1} \theta^{k-1}),\]
where $a_0, \ldots, a_{k-1}$ are integers between $-w 2^n$ and $w 2^n$. The value of $\ell_w(x)$ is taken to be $\infty$ if $x$ cannot be expressed in this form for any $n$. 

Then $\{\ell_w\}$ is a landmark system on $[0,1]$ with $r = 2^{-1}; \gamma = \sigma = k$.
\end{myex}
\begin{proof}
The monotonicity property follows because the minimum is taken over a larger set if $w$ increases.

To prove the additive property, we want to consider the sum
\[2^{-n_1} (a_0 + a_1 \theta + \cdots + a_{k-1} \theta^{k-1}) + 2^{-n_2} (b_0 + b_1 \theta + \cdots + b_{k-1} \theta^{k-1}),\]
where $|a_j| \leq w_1 2^{n_1}$ for all $j$ and $|b_j| \leq w_2 2^{n_2}$ for all $j$. Without loss of generality, we will assume $n_2 \geq n_1$. Then the sum can be rewritten as
\[2^{-n_2} ((2^{n_2 - n_1} a_0 + b_0) + (2^{n_2 - n_1}a_1 + b_1) \theta + \cdots + (2^{n_2 - n_1} a_k + b_k) \theta^k).\]

Here, each $2^{n_2 - n_1} a_j \leq w_1 2^{n_2}$, so the additivity property holds. Here, $\phi_1(w_1, w_2) = w_1 + w_2$ and there is no error term.

To prove the multiplicative property, we want to consider the product
\[2^{-n_1} (a_0 + a_1 \theta + \cdots + a_{k-1} \theta^{k-1}) \cdot 2^{-n_2} (b_0 + b_1 \theta + \cdots + b_k \theta^{k-1}),\]
where each $|a_j| \leq w_1 2^{n_1}$ and each $|b_j| \leq w_2 2^{n_2}$. The product is
\[2^{-n_1 - n_2} \sum_{j_1=0}^{k-1} \sum_{j_2 = 0}^{k-1} a_{j_1} b_{j_2} \theta^{j_1 + j_2}.\]
For $j_1 + j_2 \leq k-1$, we do not need to re-write the term; for $j_1 + j_2 \geq k$, we have that $\theta^{j_1 + j_2}$ can be expressed as a polynomial of degree at most $k-1$ in $\theta$ with integer coefficients. Therefore, the sum reduces to
\[2^{-n_1 - n_2} (c_0 + c_1 \theta + \cdots + c_k \theta^k)\]
where each of $c_0, \ldots c_k$ is a sum of a bounded number of integers (say, at most $T$) that are bounded above by $w_1 w_2 2^{n_1 + n_2}$. Thus the multiplicative property holds with $\phi_2(w_1, w_2) = Tw_1 w_2$ and no error term.

An elementary theorem (implicit in the proof of Theorem III of Chapter 5 of \cite{Cassels57}) states that there are no $k$-tuples $(a_0, \ldots, a_{k-1})$ such that $|a_0|, \ldots, |a_{k-1}| \leq 2^{n+1} w$ and 
\[|a_0 + a_1 \theta + \cdots + a_{k-1} \theta^{k-1}| \leq C (w 2^n)^{-(k-1)}.\]

This implies the separation condition.

It remains to verify the ubiquity condition. This will follow from the separation condition given above together with a transference principle \cite[Chapter 5, Theorem VI, corollary]{Cassels57}:

\begin{mythm}[Transference Principle]\label{Transference Principle}
If $\mathbf{z}$ is a $k-1$-dimensional vector such that $\mathbf{u} \cdot \mathbf{z}  - y > C_1 X^{-(k-1)}$ for all integer vectors $\mathbf{u}$ satisfying $|\mathbf{u}|_{\infty} \leq X$ and all integers $y$, then there exist constants $C_2, C_3$ such that for any real number $x$ such that $0 < x < 1$, there exists a vector $\mathbf{a}$ with $|\mathbf{a}|_{\infty} < C_2 X$ and an integer $a_0$ such that $|\mathbf{a} \cdot \mathbf{z} + a_0 - x| \leq C_3 X^{-(k-1)}$.
\end{mythm}

This transference principle, applied with $X = w 2^n$ and $\mathbf{z} = (\theta, \theta^2, \ldots, \theta^{k-1})$, immediately implies that, for any real number $x$, there exist $(a_1, \ldots, a_{k-1})$ such that the point $a_1 \theta + \cdots + a_{k-1} \theta^{k-1} - 2^nx$ is within $C w 2^{-(k-1)n}$ of an integer $-a_0$, where $|a_1|, \ldots, |a_{k-1}| < C w 2^n$ for some appropriate constant $C$ depending on $\theta$. This, of course, implies that $|a_0|$ is itself at most $C' w 2^n$, where $C'$ depends on $\theta$ but not on $n$ or on $w$. Thus, there exist $a_0, \ldots, a_{k-1}$ with $a_0, \ldots, a_{k-1} \leq C' w 2^n$ such that $a_0 + a_1 \theta + \ldots + a_{k-1} \theta^{k-1}$ is within $C w 2^{-n(k-1)}$ of $2^nx$. Dividing by $2^n$ gives the result.
\end{proof}
It is not clear if there is any way to find an appropriate landmark system for polynomials with transcendental coefficients. Nonetheless, for polynomials with coefficients well-approximated by rational numbers, we at least have access to a weak approximate landmark pair.
\begin{myex}\label{RationalApproximateExample}
Let $p(x_1, \ldots, x_v)$ be a polynomial of degree $d$ such that the coefficients of $p$ are simultaneously well-approximable to degree $\tau$; that is, $|c\mathbf{a} - \mathbf{y}| \leq C c^{-\tau}$ has infinitely many solutions for positive integers $c$ and integer vectors $\mathbf{y}$ where $\mathbf{a}$ is the vector of coefficients of $p$, and where $C$ is an appropriate constant. Let $\alpha > d/\tau$, and let $J$ be the set of values $j$ such that there exists an integer $c_j$ such that $2^{\alpha (j-1)} \leq c_j < 2^{\alpha j}$ and such that $|c_j \mathbf{a} - \mathbf{y}| \leq C 2^{\tau (\alpha + 1)j}$. Select such a $c_j$ for every $j \in J$ and define $\ell_1(x)$ to be the minimal value of $j \in J$ for which
\[x = \frac{a}{2^j}\]
for some integer $a$ between $-2^j$ and $2^{j}$ if such a $j$ exists, and $\infty$ otherwise. Define $\ell_2(x)$ to be the minimal value of $\lceil(d + \alpha)j\rceil$, where $j \in J$ is such that
\[x = \frac{a}{c_j2^{dj}}\]
for some integer $a$ if such a $j$ exists, and $\infty$ otherwise. Then $(\ell_1, \ell_2)$ is a weak-approximate landmark pair satisfying $r = \frac{1}{2}, \gamma = 1, \sigma = 1$ of degree $d + \alpha$. We emphasize the loss in the degree: although $p$ is a polynomial of degree $d$, $(\ell_1, \ell_2)$ is only a weak approximate landmark pair of degree $d + \alpha$.
\end{myex}
\begin{proof}
Let $j \in J$ and let $c_j$ be as described. Clearly multiples of $2^{-j}$ are separated by $2^{-j}$, multiples of $c_j^{-1} 2^{-d j}$ are separated by $2^{-(d + \alpha)j}$, and therefore the separation and ubiquity conditions are satisfied for $\gamma = \sigma = 1$, as described in Example \ref{RationalExample}. The factor of $d + \alpha$ was introduced in order to allow us to choose $\sigma = \gamma = 1$ for this example. Then $c_jp$ is a polynomial whose coefficients lie within $C 2^{-\alpha j \tau}$ of integers. Therefore, if we plug in numbers of the form $\frac{a}{2^j}$ into $p$, where $|a| \leq 2^j$ is an integer, we get that $c_jp(x_1, \ldots, x_v)$ is within $C' 2^{-\alpha j\tau}$ of an integer multiple of $2^{-dj}$ for an appropriate constant $C'$. Dividing by $c_j$, we get that $p(x_1, \ldots, x_v)$ is within a $C''2^{-\alpha j(\tau + 1)}$ of an integer multiple of $2^{-dj} c_j^{-1}$. The condition $\alpha > d/\tau$ implies that $\alpha j (\tau + 1) = (\alpha \tau + \alpha)j > (d + \alpha)j$, completing the proof.
\end{proof}
For similar reasons, it is possible to construct weak approximate landmark pairs for polynomials whose coefficients are simultaneously well-approximated by algebraic numbers.
\begin{myex}\label{AlgebraicApproximateExample}
Let $p(x_1, \ldots, x_v)$ be a polynomial of degree $d$ such that there exists an algebraic integer $\theta$ of degree $k$ and a real number $\tau > k - 1$ such that, for infinitely many (rational) integer vectors $(c_0, \ldots, c_{k-1})$, the polynomial $(c_0 + c_1 \theta + c_{k-1} \theta^{k-1}) p(x_1, \ldots, x_v)$ has coefficients of the form $b_0 + b_1 \theta + \cdots + b_{k-1} \theta^{k-1} + \delta$, where the coefficients $b_0, \ldots, b_{k-1}$ are integers, and $|\delta| \leq \max(|c_0|, |c_1| \ldots, |c_{k-1}|)^{-\tau}$. Let $\alpha > \frac{dk}{\tau + 1 - k}$. Let $J$ be the set of values $j$ such that there exist $c_0, \ldots, c_{k-1}$ that satisfy $2^{\alpha(j-1)} \leq \max_i |c_i| < 2^{\alpha j}$, and such that each coefficient of $(c_0 + \cdots + c_{k-1} \theta^{k-1})p(x_1, \ldots, x_v)$ is of the form described in the previous sentence. Select such a vector $\mathbf{c}^{(j)}$ for every $j \in J$, and define $\ell_1(x)$ to be the minimal value of $j \in J$ for which
\[x = \frac{1}{2^j} (a_0 + a_1 \theta + \cdots + a_{k-1} \theta^{k-1})\]
where each $a_i$ is an integer and $\max_i |a_i| \leq 2^j$ if such a $j$ exists, and $\infty$ otherwise. Define $\ell_2(x)$ to be $\lceil (d + \alpha) j \rceil$, where $j \in J$ is the minimal value such that
\[x = \frac{1}{(c_0^{(j)} + c_1^{(j)} \theta + \cdots + c_{k-1}^{(j)} \theta^{k-1}) 2^{dj}} (a_0 + \cdots + a_{k-1}  \theta^{k-1})\]
where each $a_i$ is an integer bounded above by $W 2^{(d + \alpha)j}$ in absolute value for some sufficiently large $W$ depending on $p$, $\theta$, and $k$, and $\ell_2(x) = \infty$ otherwise. Then $(\ell_1, \ell_2)$ is a weak approximate landmark pair satisfying $\gamma = k$, $\sigma = k$, and $r = \frac{1}{2}$ of degree $d + \alpha$.  
\end{myex} 
\begin{proof}
The separation and ubiquity conditions have already been verified for $\ell_1$ and $\ell_2$ with $\gamma = \sigma = k$ in Example \ref{AlgebraicExample}. We need only verify the function property.

 Let $p(x_1, \ldots, x_v)$ be a polynomial. Suppose there exists $c_0 + \cdots + c_{k-1} \theta^{k-1}$ such that $2^{\alpha (j-1)} \leq \max_i(|c_i|) < 2^{\alpha j}$ and $(c_0 + c_1 \theta + \cdots + c_{k-1} \theta^{k-1}) p(x_1, \ldots, x_v)$ has coefficients of the form $a_0 + a_1 \theta + \cdots + a_{k-1} \theta^{k-1} + \delta$, where $|\delta| \leq 2^{- \tau \alpha (j - 1)}$; this is possible whenever $j \in J$ by assumption. For each variable $x_1, \ldots, x_v$, we plug in a number of the form
\[2^{-j} (b_0 + b_1 \theta + \cdots + b_{k-1} \theta^{k-1}),\]
where $\max_i (|b_i|) \leq 2^j$, into $(c_0 + c_1 \theta + \cdots + c_{k-1} \theta^{k-1}) p(x_1, \ldots, x_v)$. Each such number satisfies $\ell_1(x) \leq j$, and each number $x$ such that $\ell_1(x) = j$ is of this form. Then, the output is within $W2^{-\tau \alpha j}$ of a number of the form $2^{-dj} (d_0 + d_1 \theta + \cdots + d_{k-1} \theta^{k-1})$ for some number $W$ that depends only on $\theta, k$ and on the polynomial $p$. Here $d_0, \ldots, d_{k-1}$ are integers that are bounded above in absolute value by $W 2^{(d + \alpha)j}$ where, again, $W$ is some constant depending only on $\theta, k,$ and $p$. We then divide by $c_0 + c_1 \theta + \cdots + c_{k-1} \theta^{k-1}$, and conclude that, for $x_1, x_2, \ldots, x_v$ with $\max(\ell_1(x_1), \ldots, \ell_1(x_v)) = j$, we have that $p(x_1, \ldots, x_v)$ is in a $2^{-(\tau + 1) \alpha j}$-neighbourhood of a number of the form 
\[\frac{2^{-d j}}{c_0 + c_1 \theta + \cdots + c_{k-1} \theta^{k-1}} (d_0 + d_1 \theta + \cdots + d_{k-1} \theta^{k-1})\]
where the $d_i$ are integers with absolute value bounded above by $W 2^{(d + \alpha)j}$. This is to say that $p(x_1, \ldots, x_v)$ is within a $2^{-(\tau + 1) \alpha j}$ neighbourhood of a point such that $\ell_2(j) \leq \lceil (d + \alpha)j \rceil$. The choice of $\alpha$ guarantees that $(\tau + 1)\alpha > (d + \alpha)k$, and because $\gamma = k$ it follows that $(\ell_1, \ell_2)$ is a weak approximate landmark pair.
\end{proof}
\section{Examples of Non-Archimedean Landmarks}
Landmark systems also arise in the setting of non-archimedean local fields, such as the $p$-adic numbers and the field of formal Laurent series over a finite field. While landmark systems for the $p$-adic numbers and for function fields are fairly simple to construct, the construction of landmark systems on other non-archimedean local fields (i.e., simple algebraic extensions of the $p$-adic numbers) is more involved. We begin by providing a landmark system for function fields, which is the simplest case.

\subsection{Introduction to Non-Archimedean Local Fields}
Before describing landmark systems for non-archimedean local fields, we briefly discuss the theory of such fields. 

A \textbf{discrete valuation ring} $R$ is a principal ideal domain with a unique prime ideal \cite{Serre79}. Because $R$ is a principal ideal domain, the prime ideal of $R$ is generated by a single element of $R$; such elements are called \textbf{uniformizers} or \textbf{uniformizing elements} of $R$. Let $t$ be a uniformizing element of $R$. Because $tR$ is the only prime ideal of $R$, it follows that $tR$ is not properly contained in any other prime ideals of $R$; therefore, $tR$ is a maximal ideal. It follows that the quotient $R/tR$ is a field. This field $R/tR$ is called the \textbf{residue class field} of $R$. We will exclusively consider the situation in which $R/tR$ is a finite field $\mathbb{F}_q$.

Suppose $q = p^f$ for some $f$. Then $\mathbb{F}_q$ has characteristic $p$, and $p \cdot 1 = \underbrace{1 + 1 + \cdots + 1}_{\text{$p$ times}}$ belongs to the ideal $tR$. If $p \cdot 1 = 0$, then the ring $R$ has characteristic $p$; otherwise, $R$ has characteristic zero.

Let $S$ be a family of representatives of the additive cosets of $tR$ in $R$ with the property that $0 \in S$. Every element $x$ of $R$ can be expressed uniquely in the form
\begin{equation}\label{DiscreteValuationRingExpansion}
x = \sum_{j=0}^{\infty} x_j t^j
\end{equation}
where $x_j$ runs over $S$. If each infinite sum of this form corresponds to an element $x \in R$, then the discrete valuation ring $R$ is called \textbf{complete}. 

For the rest of this section, we will assume $R$ is a complete discrete valuation ring. Let $x \in R$ and write $x$ as in (\ref{DiscreteValuationRingExpansion}). We define the absolute value $|x|$ of $x$ to be $0$ if $x = 0$, and $q^{-j}$ if $x_j \neq 0$ and $x_k = 0$ for all $k < j$. With respect to this absolute value, $R$ forms a complete metric space. This absolute value is discrete (this is the origin of the term \emph{discrete valuation ring}), taking only values $\{q^{-j} : j \in \mathbb{Z}\}$ and zero. The closed balls of radius $q^{-j}$ in the metric induced by this absolute value are disjoint. This absolute value respects multiplication: $|xy| = |x| |y|$. Furthermore, the absolute value satisfies the ultrametric inequality
\begin{equation}\label{UltrametricInequality}
|x + y| \leq \max(|x|, |y|).
\end{equation}
We will take a few moments to consider the importance of inequality (\ref{UltrametricInequality}). Consider the closed ball of radius $r = q^{-j}$ centered at $x$. Let $y$ be any point in $R$ such that $|x - y| \leq r$, and let $z \in R$ be such that $|y - z| \leq r$. Then we have
\[|x - z| = |(x-y)+ (y-z)| \leq \max(|x - y|, |y - z|) \leq r.\]
So the closed ball of radius $r$ centered at $x$ is precisely the same ball as the closed ball of radius $r$ centered at $y$. This implies that if two closed balls of radius $r$ intersect, then they must be equal.

The discrete nature of the absolute value also has some profound implications for the topology on $R$. For example, consider the family of closed balls of radius $q^{-j}$ contained in a closed ball of radius $q^{-(j-1)}$ centered at $x$. If $|x - y| = q^{-(j-1)}$ exactly, then $x$ and $y$ lie in the same coset of $t^{j-1}R$ but not in the same coset of $t^jR$. Since there are $q$ cosets of $t^{j}R$ contained in each coset of $t^{j-1}R$, it follows that there are precisely $q$ closed balls of radius $q^{-j}$ contained in each closed ball of radius $q^{-(j-1)}$. We can also conclude that if two open balls of radius $q^{-j}$ differ, then they are separated by a distance of at least $q^{-j}$. 

In the same spirit as for $R$, we define on norm on $R^{n}$ by 
\[\norm{(x^{(1)}, \ldots, x^{(n)})} = \max(|x^{(1)}|, \ldots, |x^{(n)}|)\]. This norm also satisfies the ultrametric property under addition, and therefore also has the property that two distinct open balls of the same radius $q^{-j}$ are separated by at least $q^{-j}$, and has the further property that each ball of radius $q^{-(j-1)}$ contains exactly $q^{n}$ balls of radius $q^{-j}$.

We now describe the Haar probability measure $dx$ on $R$: The closed ball $\overline{B(0,1)} = R$ is assigned a measure of $1$, and any closed ball of radius $q^{-j}$ is assigned a measure of $q^{-j}$. With respect to this Haar measure, any coset of $t^j R$ has measure $q^{-j}$. We will also write $dx$ for the Haar measure on $R^n$, which is the $n$-fold product of the Haar measure on $R$.

Given a complete discrete valuation ring $R$, we let $K$ be the field of fractions over $R$. Each nonzero element of $K$ is of the form
\begin{equation}\label{LocalFieldExpansion}
x = \sum_{j = M}^{\infty} x_j t^j
\end{equation}
for some possibly negative integer $M$ with $x_M \neq 0$. The field $K$ is called a \textbf{non-archimedean local field}. We extend the absolute value on $R$ to all of $K$ by defining $|x| = q^{-M}$, where $M$ is as in (\ref{LocalFieldExpansion}). We extend the Haar probability measure on $R$ to a $\sigma$-finite Haar measure on $K$ by defining the measure of a closed ball of radius $q^j$ to be $q^j$, and extend the Haar measure on $R^n$ to a $\sigma$-finite Haar measure on $K^n$ that assigns a measure of $q^{jn}$ to a closed ball of radius $q^j$.

Note that $R$ can be recovered from $K$ algebraically as the ring of integers of $K$, and topologically as the closed unit ball of $K$.

We will present two examples of non-archimedean local fields. The first such example will be the field $\mathbb{Q}_p$, known as the $p$-adic numbers. Each $x \in \mathbb{N}_0$ has a finite base-$p$ expansion
\[\sum_{j=0}^{\infty} x_j p^j\]
where only finitely many $x_j$ are nonzero. We define $|x|_p$ to be $p^{-j}$, where $x_j$ is the lowest-degree nonzero term in the expansion. If we take the completion of $\mathbb{N}_0$ with respect to this absolute value, we get the ring of elements of the form 
\[\sum_{j=0}^{\infty} x_j p^j.\]
This ring is called the \textbf{ring of $p$-adic integers}, denoted $\mathbb{Z}_p$. The $p$-adic integers are a discrete valuation ring with prime ideal $p \mathbb{Z}_p$. Any element of $\mathbb{Z}_p$ with absolute value equal to $1$ has a multiplicative inverse in $\mathbb{Z}_p$. Therefore, $\mathbb{Z}_p$ contains every rational number $\frac{r}{q}$ whose denominator $q$ is relatively prime to $p$. This is a compact abelian group under addition.

The field of fractions of $\mathbb{Z}_p$ is denoted $\mathbb{Q}_p$ and is known as the \textbf{field of $p$-adic numbers}. As an additive group, $\mathbb{Q}_p$ is locally compact.

A second example of a non-archimedean local field is the field $\mathbb{F}_q((t))$ of formal Laurent series over the finite field $\mathbb{F}_q$. Such fields are sometimes known as \textbf{function fields}. The ring of integers $\mathbb{F}_q[[t]]$ consists of formal power series over $\mathbb{F}_q$. Unlike the case for $\mathbb{Q}_p$, the field $\mathbb{F}_q((t))$ has finite characteristic $p$ where $q = p^f$. 

These two examples are central to the theory of non-archimedean local fields. According to Theorem 5 of Section 1.3 and Theorem 8 of Section 1.8 of Andr\'e Weil's book \cite{Weil95}, every non-archimedean local field is either isomorphic to $\mathbb{F}_q((t))$ or to a finite extension of $\mathbb{Q}_p$. Finite extensions of $\mathbb{Q}_p$ will be discussed in detail.

Let $L/\mathbb{Q}_p$ be a finite extension. $L$ must have residue field isomorphic to $\mathbb{F}_q$ for some $q = p^f$. The value $f$ is known as the \textbf{inertia degree} of the extension $L/\mathbb{Q}_p$. Furthermore, the absolute value on $\mathbb{Q}_p$, which will be denoted $| \cdot |_p$, has an extension to $L.$ Let $t$ be any uniformizing element of $L$. Then $|t|_p$ is equal to $p^{1/e}$ for some integer $e$. This integer $e$ is called the \textbf{ramification index} of the extension $L/\mathbb{Q}_p$. The degree of the extension $L/\mathbb{Q}_p$ is exactly $e \cdot f$. An extension $L/\mathbb{Q}_p$ is said to be \textbf{unramified} if $e = 1$, and is said to be \textbf{totally ramified} if $f = 1$. An extension $L/\mathbb{Q}_p$ always has a (not necessarily unique) subfield $K$ such that the field extension $K/\mathbb{Q}_p$ is unramified; this is called a \textbf{maximal unramified subextension} of the extension $L/\mathbb{Q}_p$.

\subsection{Landmark Systems for Non-Archimedean Local Fields}

Let $K = \mathbb{F}_q((t))$ be a function field; that is, a local field of finite characteristic, and let $R$ be the ring of integers of $K$. If $x \in R$, we can write $x$ in the form 
\begin{equation}\label{PowerSeriesExpansion}
x = \sum_{j=0}^{\infty} x_j t^j,
\end{equation}
where $x_j \in \mathbb{F}_q$ for all $j$, and $t$ is a formal variable. Addition in $R$ and multiplication in $R$ are defined in the usual way: The sum $x + y$ is defined by
\[\sum_{j=0}^{\infty} (x_j + y_j) t^j\]
and the product of $x$ and $y$ is 
\[\sum_{j=0}^{\infty} \left(\sum_{k_1 + k_2 = j} x_{k_1} y_{k_2} \right) t^j.\]
Each of the sums $\sum_{k_1 + k_2 = j} x_{k_1} y_{k_2}$ is finite and therefore is defined in $\mathbb{F}_q$.

On $\mathbb{F}_q((t))$, the absolute value of 
\[x = \sum_{j=M}^{\infty} x_j t^j\]
is equal to $q^{-M}$ if $x_M$ is nonzero. The ring $\mathbb{F}_q[[t]]$ consists of points for which the absolute value is bounded above by $1$.

The dense subring $\mathbb{F}_q[t]$ of $\mathbb{F}_q[[t]]$ consists of values $x$ with finite expansion; that is, polynomials in $t$. This dense subring gives rise to a landmark system for the compact set $\mathbb{F}_q[[t]]$.
\begin{myex}[Landmark system for $\mathbb{F}_q((t))$]\label{LandmarkSystemFunctionField}

We define a landmark function that is finite on $\mathbb{F}_q[t]$ and infinite on the rest of $\mathbb{F}_q((t))$. For $x \in \mathbb{F}_q[t]$, define $\ell(x)$ to be the minimal $j$ such that $x$ can be written in the form
\[x_0 + x_1 t + \cdots + x_j t^j\]
and $0$ if $x=0$. That is, $j$ is simply the degree of the polynomial $x$ if $x$ is a nonzero polynomial, and $\infty$ if $x$ is not a polynomial.

Then $\ell_w(x) = \ell(x)$ is a landmark system on the compact set $\Omega = \mathbb{F}_q[[t]]$ with $\gamma = \sigma = 1$. 
\end{myex}
\begin{proof}
There is nothing to check regarding the monotonicity property because $\ell_w$ is independent of $w$.

The additive and multiplicative properties are clearly satisfied because of the usual properties of the degree of a polynomial.

The separation condition is satisfied: the polynomial $x_0 + \cdots + x_j t^j$ is the only polynomial satisfying $\ell(x) \leq j$ in the open ball of radius $q^{-j}$ centered at $x_0 + \cdots + x_j t^j$.

The ubiquity condition is satisfied: the open ball of radius $q^{-s}$ centered at $x_0 + \cdots + x_s t^s$ will in fact contain exactly $q^{j - s}$ points such that $\ell(x) \leq j$.
\end{proof}

The $p$-adic integers $\mathbb{Z}_p$ consist of numbers of the form
\[x = \sum_{j=0}^{\infty} x_j p^j\]
where each $x_j \in \{0, 1, \ldots, p-1\}$. The $p$-adic numbers $\mathbb{Q}_p$ are defined similarly, but the sum is allowed to contain a finite number of terms for which $p$ has a negative power.

If $x \in \mathbb{Z}_p$ has a finite expansion, we can view $x \in \mathbb{Z}$ as a nonnegative integer in base $p$. We define a landmark function in analogy to the function field case.

\begin{myex}\label{pAdicExample}
For $x \in \mathbb{Z} \setminus \{0\}$, define
\[\ell(x_0 + \cdots + x_j p^j) = j\]
and
\[\ell(-(x_0 + \cdots + x_j p^j)) = j\]
if $x_j \neq 0$, define $\ell(0) = 0$, and take $\ell(x) = \infty$ for $x \in \mathbb{Z}_p \setminus \mathbb{Z}$.
Then $\ell_w = \ell$ is a landmark system for the compact set $\mathbb{Z}_p$ with $\gamma = \sigma = 1$.
\end{myex}

We point out that this is the first example in this paper for which there is a need for the error terms in the additive and multiplicative properties for landmark systems. For this example, the error terms can in fact taken to be bounded above by the constant function equal to $1$ everywhere.
\begin{proof}
The monotonicity property is trivial because $\ell_w$ does not depend on $w$.

The additive property follows for nonnegative integers $x$ and $y$ because
\[\sum_{i=0}^j x_i p^i + \sum_{i=0}^j y_i p^i = \sum_{i=0}^j (x_i + y_i)p^i\]
and each $x_i + y_i$ is bounded above by $2p - 2$, so $\ell(x + y) \leq \max(\ell(x), \ell(y)) + 1$. Another way of saying this is that we need only carry one digit when adding elements of $\mathbb{Z}$ in base $p$. Of course the same holds if one or both of $x$ or $y$ is negative.

The multiplicative property follows because of basic properties of arithmetic in $\mathbb{Z}$: if $x$ and $y$ are integers with absolute value strictly less than $p^{j_1 + 1}$ and $p^{j_2 + 1}$  (which is equivalent to saying that $\ell(x) \leq j_1$ and $\ell(y) \leq j_2$), then $xy$ is bounded above by $p^{j_1 + j_2 +2}$ and therefore $\ell(xy) \leq \ell(x) + \ell(y) + 1$.

The separation property holds because the ball of radius $p^{-j}$ centered at zero does not contain any other points $x$ such that $\ell(x) \leq j$: The condition that $x$ is within $p^{-j}$ of zero is equivalent to the statement that $x$ is congruent to $0$ modulo $p^{j+1}$, and the statement that $\ell(x) \leq j$ means that $x$ is an integer between $-(p^{j+1} - 1)$ and $p^{j + 1} - 1$,  and the only integer satisfying all of these conditions is $0$.

The ubiquity property holds because a $p$-adic open ball of radius $p^{-s}$ contains exactly $2p^{j - s}$ elements $y$ such that $\ell(y) \leq j$ (or $2p^{j-s} - 1$ if the $p$-adic open ball happens to contain zero).
\end{proof}

Both of these examples (as well as Example \ref{RationalExample}) rely on fairly simple algebra. For finite extensions of $\mathbb{Q}_p$, constructing a landmark system requires some nontrivial algebraic facts. 

We will first consider the case of unramified extensions of $\mathbb{Q}_p$. Let $K/\mathbb{Q}_p$ be an unramified extension of $\mathbb{Q}_p$. $K$ is formed by enlarging the residue field of $\mathbb{Q}_p$. Let $\mathbb{F}_{p^f}$ be the residue field of $K$ and let $R$ be the ring of integers of $K$. Then the field $R/pR$ is isomorphic to $\mathbb{F}_{p^f}$. We will normalize the absolute value on $K$ so that $|p|_K = p^{-f}$; this is necessary in order to guarantee that the Hausdorff dimension of $K$ is equal to $1$. Select $\alpha \in R$ such that $\alpha$ (mod $pR$) generates the multiplicative group $\mathbb{F}_{p^f}^*$ consisting of the nonzero elements of $\mathbb{F}_{p^f}$. Then $1, \alpha, \alpha^2, \ldots, \alpha^{f-1}$ are linearly independent over $\mathbb{Q}_p$. Then $\alpha$ satisfies the relation $b(\alpha) \equiv 0$ (mod $pR$) where $b(\alpha)$ is the $p^f$th cyclotomic polynomial. Furthermore, the derivative $b'(\alpha)$ is seen to have absolute value $1$, because the cyclotomic polynomial $b$ on $\mathbb{F}_{p^f}$ does not have multiple roots. Therefore, we can apply Hensel's lemma to conclude that there is an element $t \in R$ satisfying $t \equiv \alpha (\text{mod } pR)$ such that $b(t) = 0$. We know that $\{1, t, \ldots, t^{f-1}\}$ must be a basis for the extension $K/\mathbb{Q}_p$, since reducing modulo $pR$ gives a basis for $\mathbb{F}_{p^f}$ (this follows from the choice of $\alpha$). Each coefficient of $b$ is an integer and can therefore be viewed as an element of $\mathbb{Z}_p$ such that $\ell_{\mathbb{Z}_p}$ is finite. Therefore, each power $\{t^j : 0 \leq j \leq 2f - 1\}$ can be written as a $\mathbb{Z}_p$-linear combination of $1, t, \ldots, t^{f-1}$, where each coefficient maps to a finite number under $\ell_{\mathbb{Z}_p}$. 

Because $\{1, t, t^2, \ldots, t^{f-1}\}$ forms a basis of the free module $R/\mathbb{Z}_p$, we can write each element $x \in R$ in the form
\[x^{(0)} + t x^{(1)} + \cdots t^{f-1} x^{(f-1)}\]
where $x^{(k)} \in \mathbb{Z}_p$ for all $k$.

\begin{myex}[Landmark Systems for Unramified Extensions of $\mathbb{Q}_p$]\label{UnramifiedExample}
Suppose $x \in R$ can be written in the form
\[x^{(0)} + x^{(1)} + \cdots + t^{f-1} x^{(f-1)}\]
where $t$ is as constructed above. Then we define the landmark function $\ell$ by
\[\ell(x) = \max(\ell(x^{(0)}), \ldots, \ell(x^{(f-1)})).\]
Take $\ell_w(x) = \ell(x)$ for all $w$. Then $\ell_w$ is a landmark system with $\sigma = \gamma = 1$; $r = q^{-1} = p^{-f}$.
\end{myex}
\begin{proof}
The monotonicity property is trivial because $\ell_w$ does not depend on $w$.

The additive property holds because it holds in each component.

Let $x$ and $y$ satisfy $\ell(x) \leq j_1$ and $\ell(y) \leq j_2$. Then

\[xy = \sum_{i=0}^{f-1} \sum_{k=0}^{f-1} x_i y_k t^{i + k} \]

We observed above that $t^{i + k}$ can be can be written as a linear combination of $1, \ldots, t^{f-1}$ where each coefficient has finite value under $\ell_{\mathbb{Z}_p}$, so $xy$ can be written as a sum of a bounded number of terms that map to no more than $j_1 + j_2 + C'$ under $\ell$ for some constant $C'$. Applying the additive property, we conclude that $\ell(xy) \leq j_1 + j_2 + C$ for some constant $C$.

To see the separation property, notice that, for the ball of radius $q^{-j}$ containing $0$, $0$ is the only element of this ball with $\ell(x) \leq j$: writing $x = x^{(0)} + x^{(1)} t + \cdots + x^{(f-1)} t^{f-1}$, each component $x^{(i)}$ must necessarily be congruent to $0$ modulo $p^{j+1}$, and, when viewed as an integer, must be between $-(p^{j+1} - 1)$ and $p^{j+1} -1$. Therefore each component is zero. 

The ubiquity property is a consequence of the $p$-adic version: consider a ball of radius $q^{-s}$ centered at $x_0 + x_1 t + \cdots + x_{f-1} t^{f-1}$. This is an $f$-fold Cartesian product of $p$-adic balls, and therefore contains at least $(2p^{j-s} - 1)^{f} \approx_f p^{f(j-s)} = q^{j-s}$ points such that $\ell(x)$ is at most $j$.
\end{proof}
We will now extend this argument to arbitrary finite extensions $K/\mathbb{Q}_p$. Let $R_K$ be the ring of integers of $K$, and let $I_k$ be the unique prime ideal of $R_K$. Then $R_K/I_K$ is a field, and is isomorphic to $\mathbb{F}_{p^f}$, where $f$ is the inertia degree of the extension $K/\mathbb{Q}_p$. We will normalize the absolute value on $K$ so that $|s|_K = q^{-1}$ for any uniformizing element $s$ of $K$; that is, for any generator of the principal ideal $I_k$. This normalization is necessary in order to guarantee that the Hausdorff dimension of $K$ is equal to one. If $K/\mathbb{Q}_p$ is a finite extension of $\mathbb{Q}_p$, then there is an intermediate field extension $L$, a maximal unramified subextension of $K/\mathbb{Q}_p$, such that $L/\mathbb{Q}_p$ is an unramified extension and $K/L$ is a totally ramified extension. In particular, this means that the residue field of $L$ is $\mathbb{F}_{p^f}$, where $f$ is, as before, the inertia degree of $K/\mathbb{Q}_p$.

We will need to use the following algebraic fact about totally ramified extensions of $L$: every totally ramified extension $K$ of $L$ is generated by a root $s$ of an Eisenstein polynomial over $R_L$, the ring of integers of $L$. A proof of this fact can be found in \cite{Robert00}, Chapter 2, Section 4.2 for the special case $L = \mathbb{Q}_p$, but the proof extends to arbitrary $L$. This root $s$ can be chosen to be a uniformizing element of $K$; that is, each element $x$ in the ring of integers of $K$ can be written in the form
\[x = \sum_{j=0}^{\infty} x_j s^j\]
where the $x_j$ lie in a complete residue system for $L$ containing zero. Let $a(x) = x^e + a_{e-1} x^{e-1} + \cdots + a_0$ be the Eisenstein polynomial with $s$ as a root. Then each of the $a_{e-1}, \ldots, a_0 \in R_L$ is divisible by $p$, with $|a_0|_L = q^{-1}$. Thus $|a_0|_K = q^{-e}$, since we normalized the absolute value so that $|s|_K = q^{-1}$.

Consider the equation $a(s) = 0$, which holds in $K$. Expanding the left side of the equation, we get $s^e + a_{e-1}s^{e-1} + \cdots + a_1s + a_0 = 0$. The derivative $a'$ is nonzero at $s$; otherwise $s$ would have degree less than $e$ over $L$, which is impossible because $a$ irreducible by Eisenstein's criterion. Suppose that the absolute value of $a'(s)$ in $K$ is equal to $q^{-\alpha}$. Reduce the equation $a(s) = 0$ modulo $s^{2 \alpha + 2} R$. Performing the reduction, we get $s^e + a_{e-1} s^{e-1} + \cdots + a_1 s + a_0 \equiv 0$ (mod $s^{2 \alpha + 2} R$). This equation continues to hold if we replace $a_0, \ldots, a_{e-1}$ by any other coefficients that are congruent to $a_0, \ldots a_{e-1}$ modulo $s^{2 \alpha + 2}R$. In particular, we can replace them with elements $b_0, \ldots, b_{e-1}$ of $L$ such that $\ell_L(b_i)$ is no more than $\left\lceil \frac{2 \alpha + 2 }{e} \right\rceil$. Let $b(x)$ be a polynomial with this replacement made. Then $b(s) \equiv 0$ (mod $s^{2\alpha + 2}R$) and $|b'(s)| = q^{\alpha}$, so by the version of Hensel's lemma appearing in Chapter 2, Section 1.5 of \cite{Robert00}, it follows that $b(x)$ has a root in an open $q^{-\alpha -1}$-neighbourhood of $s$. In particular, this root must have absolute value $q^{-1}$. Let $t$ be this root of $b(x)$. Because $t$ is a uniformizer of $K$, it follows that $1, t, \ldots, t^{e-1}$ form a vector space basis for the extension $K/L$, and we can write every element of the ring of integers of $K$ in the form
\[x = x^{(0)} + x^{(1)} t + \cdots + x^{(e-1)} t^{e-1}\]
where each of the $x^{(0)}, \ldots, x^{(e-1)}$ are in the ring of integers of $L$. This number $t$ allows us to define a landmark system.
\begin{myex}
For this choice of $t$, define $\ell_K(x)$ to be 
\[\max(\ell_L(x^{(0)}), \ldots, \ell_L(x^{(e-1)})). \] 
Then $\ell_w(x) = \ell_K(x)$ is a landmark system for $R_K$ with $\gamma = \sigma = 1$ and $r = q^{-1}$.
\end{myex}
\begin{proof}
The monotonicity property is trivial because $\ell_K$ does not depend on $w$.

The additive property is shown as follows. Suppose $\ell_K(x) = j_1$ and $\ell_K(y) = j_2$. Then, writing $x = x^{(0)} + x^{(1)}t + \cdots + x^{(e-1)} t^{e-1}$ and $y = y^{(0)} + y^{(1)}t + \cdots + y^{(e-1)}t^{e-1}$, we get
\[x + y = (x^{(0)} + y^{(0)}) + (x^{(1)} + y^{(1)})t + \cdots + (x^{(e-1)} + y^{(e-1)})t^{e-1}\]
and therefore 
\begin{IEEEeqnarray*}{rCcl}
\ell_K(x+y) & = & &\max(\ell_L(x^{(0)} + y^{(0)}), \ell_L(x^{(1)} + y^{(1)}), \ldots, \ell_L(x^{(e-1)} + y^{(e-1)})) \\
& \leq & & \max \bigl( \max(\ell_L(x^{(0)}), \ell_L(y^{(0)})) + 1, \max(\ell_L(x^{(1)}) \ell_L(y^{(1)})) + 1, \\
& &,  & \ldots, \max(\ell_L(x^{(e-1)}), \ell_L(y^{(e-1)}) ) + 1 \bigr) \\
& \leq & &  \max(\ell_K(x), \ell_K(y)) + 1
\end{IEEEeqnarray*}
which shows the additive property.

Next, we show the multiplicative property. Suppose that $\ell_K(x) = j_1$ and $\ell_K(y) = j_2$. Expanding the product, we get
\[xy = \sum_{j=0}^{e-1} \sum_{k=0}^{e-1} x^{(j)} y^{(k)} t^{j+k}\]
The number of summands depends only on the field $K$ (and in particular on the ramification index $e$ of the extension $K/L$). Furthermore, each $t^{j+k}$ can be written as an $L$-linear combination of $\{1, t, \ldots, t^{e-1}\}$ where each coefficient maps to a finite number under $\ell_L$. Since each $\ell_L(x^{(j)} y^{(k)})$ is bounded above by $\ell_L(x^{(j)}) + \ell_L(y^{(k)}) \leq \ell_K(x) + \ell_K(y)$, it follows from the additive property for $\ell_K$ that the sum is bounded above by $\ell_K(x) + \ell_K(y) + C_K$ for some appropriate constant $K$.

The separation property is shown in a similar way to the unramified case. 

The ubiquity property is a consequence of the ubiquity property for unramified extensions in exactly the same way that the ubiquity property for unramified extensions of $\mathbb{Q}_p$ follows from the ubiquity property for $\mathbb{Q}_p$.
\end{proof}
\section{Avoidance of Landmark Configurations at a Single Scale: Nondegenerate Case}
Before embarking on our proof, we will make an observation. Suppose that $K$ is either $\mathbb{R}$ or some non-archimedean local field with residue field $\mathbb{F}_q$. We will briefly consider the roles of $r, \gamma,$ and $\sigma$ in the definition of weak approximate landmark pairs. We observe that if $(\ell_1, \ell_2)$ is a weak approximate landmark pair with parameters $(r, \gamma, \sigma)$, then for any $\beta > 0$, $(\ell_1, \ell_2)$ is also a weak approximate landmark pair with parameters $(r^{\beta}, \gamma/\beta, \sigma/\beta)$. Thus we can assume without loss of generality that $r = 2^{-1}$ if $K = \mathbb{R}$, or $r = q^{-1}$ if $K$ is a non-archimedean local field with residue field $\mathbb{F}_q$. In particular, we can always assume $r^{-1}$ is an integer. We will make this assumption henceforth.

We begin by considering functions with some nonzero first-order partial derivative.
\begin{mypro}\label{NondegenerateProposition}
Let $f : K^{nv} \to K$, where $K$ is either a non-archimedean local field or $\mathbb{R}$. Let $T_1, \ldots, T_v$ be compact subsets of $K^n$, each of which is a union of essentially disjoint closed cubes of sidelength $r^s$, and let the strictly differentiable function $f(x^{(1)}, \ldots, x^{(v)})$ satisfy $\left|\frac{\partial f}{\partial x^{(k_0)}_{i_0}}\right| \geq c$ for some $i_0$ and for all $(x^{(1)}, \ldots, x^{(v)}) \in T_1 \times \cdots \times T_v$. Suppose $(\ell_1, \ell_2)$ is a weak approximate landmark pair adapted to the function $f$ on a set containing the projection of each $T_j$ onto each coordinate axis. Then there exists a small positive real number $\epsilon^*$ such that for all $0 < \epsilon < \epsilon^*$, there exists $c'(\epsilon)$ with the following property. If $j \in J$ is sufficiently large, then there exist sets $S_1 \subset T_1, \ldots, S_v \subset T_v$ such that:
\begin{enumerate}
\item There are no solutions to $f(x^{(1)}, \ldots, x^{(v)}) = 0$ with $x^{(1)} \in S_1, \ldots, x^{(v)} \in S_v$. Furthermore, $f$ satisfies the lower bound $|f(x^{(1)}, \ldots, x^{(v)})| \geq c' r^{(\gamma + \epsilon) dj}$ on $S_1 \times \cdots \times S_v$. 
\item Let $U$ be one of the cubes of sidelength $r^{s}$ that constitute $T_k$. Then $S_k \cap U$ is a union of at least $c' r^{-n (\sigma - 5 \epsilon) j}$ disjoint cubes of sidelength $c'r^{d j (\gamma + \epsilon)}.$ Furthermore, for each integer $s'$ such that $s < s' \leq \sigma j - 1000 \epsilon j$, each cube of sidelength $r^{s'}$ contained in $U$ will contain at least $c' r^{-n ((\sigma - 5 \epsilon) j - s')}$ cubes of sidelength $c' r^{dj (\gamma + \epsilon)}$. We can further guarantee that each cube of sidelength $r^{\sigma j}$ will intersect no more than one such cube.
\end{enumerate}
\end{mypro}
\begin{proof}
Throughout this argument, all constants named $c$ or $C$ will depend on $\epsilon, f$, and $s$, but the dependence will be suppressed. Several measures will be taken in order to guarantee that certain points do not lie near the boundary of a cube; these precautions are unnecessary in non-archimedean local fields and can thus be ignored for that setting.

Throughout this proof, we will define $\ell_w(\mathbf{x}) = \max_k \ell_w(x_k)$ for any vector $\mathbf{x}$ of any dimension. Consider a cube $U$ of sidelength $r^s$ contained in $T_i$. Given a point $y \in U$ satisfying $\ell_1(y) \leq j$, we will define $B_y$ to be the box of sidelength $r^{(\gamma + \epsilon)dj}$ centered at $y$. 

We partition $U$ into cubes of sidelength $r^{\lceil (\sigma - 4 \epsilon) j\rceil}$. Let $V$ be such a cube. In the Euclidean setting, we will take $V'$ to be the slightly smaller cube of sidelength $r^{\lceil (\sigma - 3 \epsilon)j\rceil}$ with the same center as $V$; for non-archimedean local fields we will simply select $V' = V$. The ubiquity condition guarantees that, if $j$ is large enough,  $V'$ will contain at least one point $y$ such that $\ell_1(y) \leq j$. For each $V$, we pick such a point $y(V)$. Let $Y(U)$ be the set $\{y(V) : V \text{ is one of the cubes forming $U$}\}$. Then the cubes $\{B_y : y \in Y(U)\}$ are disjoint (as they have sidelength $r^{(\gamma + \epsilon) dj} \ll r^{(\sigma - 3 \epsilon) j }$), and each cube of sidelength $r^{\sigma j}$ intersects only one such cube provided that $j$ is large enough with respect to $\epsilon$. Let $Y_k$ be the union of the sets $Y(U)$ over all of the cubes $U$ that constitute $T_k$.

Let $\mathbf{y} := (y^{(1)}, \ldots, y^{(v)})$ where $y^{(1)} \in Y_1, \ldots, y^{(v)} \in Y_v$. In particular, $\mathbf{y}$ satisfies $\ell_w(\mathbf{y}) \leq j$. Consider the behaviour of $f$ on the product $B_{\mathbf{y}} := B_{y^{(1)}} \times \cdots \times B_{y^{(v)}}$. Because $f$ has a bounded gradient (say, bounded by $C_1$) on $B_{\mathbf{y}}$ where the bound does not depend on $\mathbf{y}$, there exists a constant $C_1 \geq 1$ such that $B_{\mathbf{y}}$ maps into a box of side length at most $C_1 r^{(\gamma + \epsilon) d j}$.  Because $(\ell_1, \ell_2)$ is a weak approximate landmark pair for $f$, it follows that if $\epsilon^*$ is small enough, $f(y^{(1)}, \ldots, y^{(v)})$ is a within a $C_2r^{(\gamma + 2\epsilon) d j}$-neighbourhood of a point $z$ satisfying $\ell_2(z) \leq dj + o(j)$. If $j \in J$ is sufficiently large, we can guarantee both that $\ell_2(z) \leq (d + \epsilon) j$  and that $C_2 r^{(\gamma + 2 \epsilon) dj } < c/8 r^{(\gamma + \epsilon) dj}$. Furthermore, if $j$ is sufficiently large depending on $\epsilon$, the separation condition implies that either the image of $B_{\mathbf{y}}$ avoids a $C_1 r^{(\gamma + \epsilon) dj}$ neighbourhood of $0$, or that $0$ is the only point $z$ in the image of $B_{\mathbf{y}}$ satisfying $\ell_2(z) \leq (d + \epsilon) j$. In particular, this means that $|f(y)| \leq \frac{c}{8} r^{(\gamma + \epsilon) dj}$.

Let $i_0, k_0$ be such that $\left| \frac{\p f}{\p x_{i_0}^{(k_0)}} \right| \geq c$. We will define $S_k$ to be a union of cubes defined as follows. If $y \in Y_k$ for some $k \neq k_0$, let $y^* = y$. If $y \in Y_{k_0}$, we instead let $y^{*} = y + \frac{1}{2} r^{(\gamma + \epsilon) dj} e_{i_0}$, where $e_{i_0}$ is the vector with a $1$ in the $i_0$ component and zeroes elsewhere. In either case, take $S_k$ to be the union over all $y \in Y_k$ of the $B_y^*$ , where $B_y^*$ is the box centered at $y^*$ with sidelength $c^* r^{\lceil (\gamma + \epsilon)dj \rceil}$, where $c^*$ is the largest integer power of $r$ that is less than $\frac{c}{4 C_1 \sqrt{n}}$. For $\mathbf{y} = (y^{(1)}, \ldots, y^{(v)})$ define $B_{\mathbf{y}}^* := B_{y^{(1)}}^* \times \cdots \times B_{y^{(v)}}^*$.

We verify that the sets $S_k$ satisfy the conditions of the Proposition. We will begin with part 1. Suppose $\mathbf{x} = (x^{(1)}, \ldots, x^{(v)})$ where $x^{(1)} \in S_1, \ldots, x^{(v)} \in S_v$. We would like to show a lower bound on $|f(x^{(1)}, \ldots, x^{(v)})|$. First, we observe that $\mathbf{x} \in B_{\mathbf{y}}^*$ for some $\mathbf{y} \in Y_1 \times \cdots \times Y_k$. We split into two cases depending on whether $f(B_{\mathbf{y}})$ contains $0$.

If $f(B_{\mathbf{y}})$ does not contain $0$, then, because $\mathbf{x} \in B_{\mathbf{y}}^* \subset B_{\mathbf{y}}$, it follows that $|f(\mathbf{x})| \geq C_1 r^{(\gamma + \epsilon)dj}$. 

If, instead, $f(B_{\mathbf{y}})$ does contain $0$, then we make use of the choice of $\mathbf{y}^*$. Note that $\mathbf{y}^* - \mathbf{y} = \frac{1}{2} r^{(\gamma + \epsilon) dj} e_{i_0}^{(k_0)}$, where $e_{i_0}^{(k_0)} = (0, \ldots, 0, e_{i_0}, 0, \ldots, 0)$ with $e_{i_0}$ in the $k_{0}$ component, and the $n$-dimensional $0$ vector in the remaining $v-1$ components. We have, by assumption, a lower bound of $c$ on the absolute value of the derivative of $f$ in the $e_{i_0}^{(k_0)}$ direction, and a bound of $C_1$ on the gradient of $f$. Because $f$ is strictly differentiable, we have that for $j$ sufficiently large,

\begin{IEEEeqnarray*}{rCl}
f(\mathbf{x}) & \geq & f(\mathbf{y}^*) - \frac{17 C_1}{16} |\mathbf{x} - \mathbf{y}^*|  \\
 & \geq & \frac{7}{8} |f'(\mathbf{y})| |\mathbf{y} - \mathbf{y}^*| - |f(\mathbf{y})| - \frac{17 C_1}{16} |\mathbf{x} - \mathbf{y}^*|\\
 & \geq &  \frac{7c}{16} r^{(\gamma + \epsilon)dj} - \frac{c}{8} r^{(\gamma + \epsilon)dj} - \frac{17 C_1c \sqrt{n}}{64 C_1 \sqrt{n}} r^{(\gamma + \epsilon)dj} \\
& \geq & \frac{3c}{64} r^{(\gamma + \epsilon)dj}
\end{IEEEeqnarray*}
Therefore, $|f(\mathbf{x})| \geq \frac{3c}{64} r^{(\gamma + \epsilon) dj}$, as desired. This establishes conclusion 1.

We now prove conclusion 2. Let $U$ Be a constituent cube of $T_k$. The number of cubes $V$ in the decomposition above is $r^{-n (\lceil (\sigma - 4 \epsilon)j \rceil - s) }$, and each cube $V$ contains a cube $B_y*$, where $y = y(V)$.  If $\tilde U$ is an arbitrary cube of sidelength $r^{s'}$ contained in $U$, then, provided that $j$ is sufficiently large, $\tilde U$ entirely contains at least $r^{-n( \lceil (\sigma - 5 \epsilon) j \rceil -s')}$ cubes $B_y*$ as desired.
\end{proof}
\section{Construction of the Set: Avoiding General Landmark Configurations at Multiple Scales}
We now construct the set $E$ promised by the statement of Theorem \ref{LandmarkTheorem}. We adopt a queueing strategy similar to the one described in \cite{FraserPramanik18} in order to construct our set. Without loss of generality, we can assume, possibly by modifying $\gamma$ or $\sigma$ if necessary, that $r$ is as described at the beginning of Section 5.

\paragraph{Stage 0} Let $E_0 = B$ where $B$ is as defined as in the statement of \ref{LandmarkTheorem}. We can assume $B$ is a closed cube of sidelength $r^{s_0}$ for some $s_0$. Fix a sequence $\epsilon_j$ such that $\epsilon_j \to 0$ and $\epsilon_j < \frac{\sigma}{1000000}$ for all $j$. Select $L_0$ sufficiently large so that the ball $E_0$ can be partitioned into at least $v_1 + 1$ essentially disjoint cubes of sidelength $r^{L_0}$. Let $B_1^{(0)}, \ldots, B_{M_0}^{(0)}$ be an enumeration of the cubes of sidelength $r^{L_0}$ contained in $B^n$, and let $\Sigma_0$ be the family of $v_1$-tuples of distinct such cubes, ordered lexicographically and identified in the usual way with the family of injections from $\{1, \ldots, v_1\}$ into $\{1, \ldots, M_0\}$. Let $\mathcal{Q}_0$ be the queue consisting of the $4$-tuples
\[\{(1, k, \tau, 0) : 0 \leq k \leq |\alpha_1| - 1, \tau\in \Sigma_0\},\]
where the queue elements are ordered so that $(1, k, \tau, 0)$ precedes $(1, k', \tau', 0)$ whenever $\tau <  \tau'$, and $(1, k, \tau, 0)$ precedes $(1, k', \tau, 0)$ whenever $k > k'$.

\paragraph{Stage 1} At Stage 1, we will consider the first queue element $(1, k, \tau, 0)$. Let $T_i^{(1)} = B_{\tau(i)}^{(0)}$ for all $1 \leq i \leq v_1$.

Let $f = D_k f_1$. By the ordering of $\mathcal{Q}_0$, we know that $k = |\alpha_1| - 1$ and therefore $\frac{\partial f}{\partial x_{i_1}^{(j_1)}} = D_{k+1} f_1$ is nonzero for some $i_1, j_1$. Furthermore, by compactness, we know that there is a lower bound, say, $r^{A_1}$ for this derivative on $B^n$. We select a weak approximate landmark pair $(\ell_1, \ell_2)$ adapted to $f$ for which the degree is at most $d + \epsilon_1$. Now we apply Proposition \ref{NondegenerateProposition} to arrive at sets $S_1^{(1)}, \ldots, S_{v_1}^{(1)}$. Let $\epsilon_1^*$ be the minimum of $\epsilon_1$ and the value $\epsilon^*$ required to apply the proposition. We can select $N_1$ ($j$ in the Proposition) to be a large number depending on $r, N_0, n, \sigma, d, \gamma,$ and $\epsilon_{1}^*$. The exact requirements on $N_1$ will be specified later. 

Then $S_1^{(1)} \subset T_1^{(1)}, \ldots, S_{v_1}^{(1)} \subset T_{v_1}^{(1)}$ have the property that $D_k f_1$ is nonzero for $x_i \in S_i^{(1)}$, where $1 \leq i \leq v_1$. We will define a subset $E_1 \subset E_0$ in the following way. We take $E_1 \cap T_1^{(1)} = E_0 \cap S_1^{(1)}, E_1 \cap T_2^{(1)} = E_0 \cap S_2^{(1)}, \ldots, E_1 \cap T_{v_1}^{(1)} = E_0 \cap S_{v_1}^{(1)}$. We decompose each of the $r^{L_0}$-cubes not contained in $T_1^{(1)} \cup \cdots \cup T_{v_1}^{(1)}$ into $r^{L_1}$-cubes, and retain all of these subcubes that do not border $T_1^{(1)} \cup \cdots \cup T_{v_1}^{(1)}$ as part of $E_1$. This gives a subset $E_1 \subset E_0$ that can be expressed as an essentially disjoint union of cubes of side length $r^{L_1}$, where $L_1$ is such that $r^{L_1} = c'r^{(d + \epsilon_1) (\gamma + \epsilon_1^*) N_1}$. We can assume that $L_1$ is an integer by shrinking the cubes from the proposition slightly if necessary. 

Let $\mathcal{E}_1$ be the collection of cubes of side length $r^{L_1}$ whose union is $E_1$. Enumerate the cubes of $\mathcal{E}_1$ as $B_1^{(1)}, \ldots, B_{M_1}^{(1)}$. For $q = 1, 2$ define $\Sigma_1^{(q)}$ to be the collection of $v_q$-tuples of distinct such cubes, ordered lexicographically and identified in the usual way with the family of injections from $\{1, \ldots, v_q\}$ into $\{1, \ldots, M_1\}$. We assume $N_1$ has been chosen sufficiently large that $\Sigma_1^{(q)}$ is nonempty for $q = 1, 2$. We then form the queue $\mathcal{Q}_1'$ consisting of $4$-tuples of the form
\[\{(q, k,\tau, 1): 1 \leq q \leq 2; 0 \leq k \leq |\alpha_q| - 1; \tau\in \Sigma_1^{(q)} \} \]
arranged so that $(q, k,\tau, 1)$ precedes $(q', k',\tau', 1)$ if $q \leq q'$, so that $(q, k,\tau, 1)$ precedes $(q, k',\tau', 1)$ if $\tau<\tau'$, and so that $(q, k,\tau, 1)$ precedes $(q, k',\tau, 1)$ if $k > k'$. We arrive at the queue $\mathcal{Q}_1$ by appending the queue $\mathcal{Q}_1'$ to $\mathcal{Q}_0$.

\paragraph*{Stage $j$} We will now describe Stage $j$ of the construction for $j > 1$. We follow essentially the same procedure as in Stage 1. We begin with a decreasing family of sets $E_0, \ldots, E_{j-1}$. Each $E_{j'}$ is a union of cubes of sidelength $r^{L_{j'}}$, the collection of which is called $\mathcal{E}_{j'}$. The family of $v_{q'}$-tuples of distinct cubes in $\mathcal{E}_{j'}$ will be denoted $\Sigma_{j'}^{(q')}$. We have a queue $\mathcal{Q}_{j-1}$ consisting of $4$-tuples $(q', k',\tau', j')$ where we have $0 \leq j' \leq j-1$, $1 \leq q \leq j'+1$, $0 \leq k' \leq |\alpha_{q'}| - 1$, and $\tau' \in \Sigma_{j'}^{(q')}$. The set $E_{j-1}$ has the property that $D_{k'} f_{q'}(x_1, \ldots, x_{v_{q'}}) \neq 0$ for $x_1 \in B_{\tau'(1)}^{(j')} \cap E_{j-1}, \ldots, x_{v_{q'}} \in B_{\tau'(v_{q'})}^{(j')} \cap E_{j-1}$ for any $(q', k',\tau', j')$ in the first $j - 1$ elements of the queue $\mathcal{Q}_{j-1}$. 

Consider the $j$th queue element $(q, k,\tau, j_0)$, where $q \leq j_0 \ll j$. We will consider two cases: the case in which $k = |\alpha_q|-1$, and the case in which $k < |\alpha_q| -1$.
\subparagraph*{Case 1: $k = |\alpha_q| - 1$.} Let $f = D_{|\alpha_{q}| - 1} f_q$. In this case, we have that $k + 1 = |\alpha_q|$; therefore, it follows by assumption that $D_{|\alpha_q|}f_q$ is nonzero on all of $T_1 \times \cdots \times T_v$. Let $r^{A}$ be the lower bound on this partial derivative. By assumption, we have a weak approximate landmark pair $(\ell_1, \ell_2)$ for $f$ of degree less than $d + \epsilon_j$. We then apply Proposition \ref{NondegenerateProposition} with the quantity $N_j$ ($j$ in the proposition) taken to be a large number depending on $r, N_{j-1}, n, \sigma, d, \gamma$, and $\epsilon_j^*$. Here, $\epsilon_j^*$ is the minimum of $\epsilon_j$ and the value of $\epsilon^*$ occurring in the proposition. The specific requirements for the choice of $N_j$ will be specified later.
\subparagraph*{Case 2: $k < |\alpha_q| - 1$.} Let $f = D_k f_q$. If $k < |\alpha_q| - 1$, then, by the ordering of the elements of the queue $\mathcal{Q}_{j-1}$, we will have that the $j-1$st element of $\mathcal{Q}_{j-1}$ is $(q, k+1,\tau, j_0)$. Therefore, by the previous stage, we have that for $x_1 \in T_1, \ldots, x_v \in T_v$ that $D_{k+1} f_q$ is nonzero. But this implies by compactness that there exists some $A$ such that $D_{k+1} f_q$ is at least $r^{-A}$ in absolute value on all of $T_1 \times \cdots \times T_v$. Furthermore, we have assumed that $f$ has a weak approximate landmark pair of degree at most $d + \epsilon_j$. Apply Proposition \ref{NondegenerateProposition} to the sets $T_1, \ldots, T_v$ with the quantity $N_j$ ($j$ in the proposition) chosen to be a very large number depending on $r, N_{j-1}, n, \sigma, d, \gamma$, and $\epsilon_{j}^*$. Here, $\epsilon_j^*$ is the minimum of $\epsilon_j$ and the value $\epsilon^*$ required to apply the proposition. The specific requirements for the choice of $N_j$ will be specified later.

In any case, we arrive at sets $S_1^{(j)} \subset T_1^{(j)}, \ldots, S_{v_q}^{(j)} \subset T_{v_q}^{(j)}$, such that $D_k f_q$ is nonzero for $(x^{(1)}, \ldots, x^{(v_q)}) \in S_1^{(j)} \times \cdots \times S_{v_q}^{(j)}$. We define a subset $E_j \subset E_{j-1}$ in the following way. We take $E_j \cap T_1^{(j)} = E_{j-1} \cap S_1^{(j)}$, $E_j \cap T_2^{(j)} = E_{j-1} \cap S_2^{(j)}, \ldots, E_{j-1} \cap T_{v_q}^{(j)} = E_j  \cap S_{v_q}^{(j)}$. We split the cubes of sidelength $r^{L_{j-1}}$ not contained in $T_1^{(j)} \cup \cdots \cup T_v^{(j)}$ into cubes of sidelength $r^{L_{j}}$; the cubes that do not border $T_1^{(j)} \cup \cdots \cup T_{v_q}^{(j)}$ will be retained as part of $E_j$. This gives a subset $E_j \subset E_{j-1}$ that can be expressed as a disjoint union of cubes of sidelength  $r^{L_j}$, where $L_j$ is the smallest integer such that $r^{L_j} \leq c' r^{(d + \epsilon_j) (\gamma + \epsilon_j^*) N_j}$. Call the collection of such balls $\mathcal{E}_j$, and let $B_1^{(j)}, \ldots, B_{M_j}^{(j)}$ be an enumeration of the balls in $\mathcal{E}_j$. For each $1 \leq q \leq j$, we define $\Sigma_j^{(q)}$ to be the collection of $v_q$-tuples of distinct balls in $\mathcal{E}_j$. We assume that $N_j$ has been chosen sufficiently large in order to guarantee that these sets will be nonempty. We equip $\Sigma_j^{(q)}$ with the lexicographic order and identify $\Sigma_j^{(q)}$ with the set of injections from $\{1, \ldots, v_q\}$ into $\mathcal{E}_j$. Consider the queue $\mathcal{Q}_j'$ consisting of $4$-tuples $(q, k,\tau, j)$ where $1 \leq q \leq j + 1$, $0 \leq k \leq |\alpha|_q - 1$, and $\tau \in \Sigma_j^{(q)}$. We order the queue $\mathcal{Q}_j'$ in the following way: $(q, k,\tau, j)$ will precede $(q', k',\tau', j)$ if $q < q'$, $(q, k,\tau,j)$ precedes $(q, k',\tau', j)$ if $\tau<\tau'$, and $(q, k,\tau, j)$ precedes $(q, k',\tau, j)$ if $k > k'$. We append the queue $\mathcal{Q}_j'$ to $\mathcal{Q}_{j-1}$ to arrive at the queue $\mathcal{Q}_j$.

The set $E$ is given by $E = \cap_{j=1}^{\infty} E_j$.
\subsection{Hausdorff Dimension of E}
We now outline the computation of the Hausdorff dimension of the set $E$. In order to compute the Hausdorff dimension of this set, we use a version of Frostman's lemma. The goal is to construct a Borel probability measure $\mu$ supported on $E$ such that $\mu(I) \lesssim_{\epsilon} l(I)^{\frac{n \sigma}{d \gamma} - \epsilon}$ for every cube $I$ with side length $l(I)$. The existence of such a measure would imply that the Hausdorff dimension of $E$ is at least $\frac{n \sigma}{d \gamma}$.

We will now describe the construction of the measure $\mu$. $\mu$ will be obtained as a weak limit of measures $\mu_j$ supported on the sets $E_j$. We begin by defining $\mu_0$ to be the uniform probability measure on the set $E_0$. Decompose $E_0$ into subcubes of sidelength $r^{\lfloor (\sigma - 1000 \epsilon_1^*) N_1 \rfloor}$, and split the mass of $E_0$ evenly among such cubes. 

Let $J$ be such a cube. Then there are two possibilities: either $J$ is contained in some $T_i^{(1)}$ for some $1 \leq i \leq v_{q_1}$, or $T_i$ is essentially disjoint from $T_i^{(1)}$ for all $i$. If $J$ is contained in a cube of $T_i^{(1)}$, then part 2 of Proposition \ref{NondegenerateProposition} states that there are at least $c'r^{-995 \epsilon_1^* N_1}$ and at most $r^{-1000 \epsilon_1^* N_1}$ cubes of radius $r^{L_1}$ contained in $J$ that are retained as part of $E_1$. 

The second case is the case in which the cube $J$ is essentially disjoint from the sets $T_i^{(1)}$. In this case, all of the subcubes of $J$ of side length $r^{L_1}$ that do not border any of the sets $T_i^{(1)}$ are retained. If $N_1$ is chosen sufficiently large, this will be at least half of the subcubes of $J$ of sidelength $r^{L_1}$. The measure $\mu_1|_J$ is obtained by splitting the measure of $J$ evenly among each of these surviving cubes.

We continue this procedure inductively. Suppose we have a probability measure $\mu_j$ supported on $E_j$. The set $E_j$ is a union of cubes of sidelength $r^{L_j}$. We will describe the construction of the measure $\mu_{j+1}$ from $\mu_j$ as follows. We will decompose each of the $r^{L_j}$-cubes that constitute $E_j$ into a union of essentially disjoint cubes of side length $r^{\lfloor (\sigma - 1000 \epsilon_{j+1}^*) N_{j+1} \rfloor}$. Each such cube will receive the same share of the parent cube's measure.

Let $J$ be one of these cubes of sidelength $r^{\lfloor (\sigma - 1000 \epsilon_{j+1}^*) N_{j+1} \rfloor}$. There are two possibilities: either $J$ is contained in some $T_i^{(j+1)}$, or $J$ is essentially disjoint from all of the sets $T_i^{(j+1)}$, for $1 \leq i \leq v_{q_{j+1}}$. 

If $J$ is contained in some $T_i^{(j+1)}$, then $J \cap E_{j+1}$ is a union of cubes of sidelength $r^{L_{j+1}}$. We observe that the number of such cubes contained in $J$ is at at least $c' r^{-995 \epsilon_{j+1}^* N_{j+1}}$ and at most $r^{- 1000 \epsilon_{j+1}^* N_{j+1}}$. If $J$ is not contained in any of the $T_i^{(j+1)}$, then $J$ is essentially disjoint from the $T_i^{(j+1)}$, and, provided $N_{j+1}$ is large enough, at least half of the subcubes of $J$ of sidelength $r^{L_{j+1}}$ are retained. In either case, we distribute the measure of $J$ evenly among all of the surviving subcubes of sidelength $r^{L_{j+1}}$ contained in $J$. 

We claim that the measures $\mu_j$ have a weak limit $\mu$, which satisfies the Frostman condition. It is clear that the measures $\mu_j$ have a weak limit because they are defined via a mass-distribution process. We will show that this weak limit $\mu$ satisfies the Frostman condition. First, we will show that the Frostman condition with dimension $\frac{n\sigma }{d \gamma}$ is satisfied for the basic cubes in the construction.

Let $I \in \mathcal{E}_j$, and let $J \in \mathcal{E}_{j+1}$ be contained in $I$. We will consider two cases: the case in which $I \subset T_i^{(j+1)}$ for some $i$, and the case in which $I$ is essentially disjoint from the sets $T_i^{(j+1)}$. We observe that it is enough to obtain an estimate for $\mu_{j+1}(J)$, because $\mu_{j'}(J) = \mu_{j+1}(J)$ for all $j' \geq j+1$.

\paragraph*{Case 1: $I$ is contained in $T_i^{(j+1)}$ for some $i$.}
In this case the measure of $I$ is split evenly among the subcubes of sidelength $r^{\lfloor(\sigma - 1000 \epsilon_{j+1}^*) N_{j+1}\rfloor}$. Each of these subcubes will therefore have measure $\mu_j(I) r^{n(L_j - \lfloor (\sigma - 1000 \epsilon_{j+1}^*) N_{j+1} \rfloor)}.$ Each such cube will contain at least $c' r^{- \lfloor 995n \epsilon_{j+1}^* N_{j+1} \rfloor}$ cubes with the same $\mu_{j+1}$-measure as $J$. Thus the $\mu_{j+1}$-measure of $J$ is at most 
\[c'^{-1} \mu_j(I) r^{-n L_j + n \lfloor (\sigma - 1000 \epsilon_{j+1}^*) N_{j+1} \rfloor + \lfloor 995 n \epsilon_{j+1}^* N_{j+1} \rfloor}.\]
After combining terms, we get an estimate of 
\[\mu_{j+1}(J) \leq c'^{-1} \mu_j(I) r^{-n L_j + n(\sigma -5 \epsilon_{j+1}^*) N_{j+1} - n}.\]
We can choose $N_{j+1}$ sufficiently large so that $c'^{-1} \mu_j(I) r^{-nL_j - n} \leq r^{-n \epsilon_{j+1}^* N_{j+1}}$. Then we get the estimate
\[\mu_{j+1}(J) \leq r^{n(\sigma - 6 \epsilon_{j+1}^*) N_{j+1}}.\]
But $J$ has sidelength $r^{(d + \epsilon_{j+1}) (\gamma + \epsilon_{j+1}^*) N_{j+1}}$.  Thus
\[\mu_{j+1}(J) \leq \ell(J)^{\frac{n(\sigma - 6 \epsilon_{j+1}^*)}{(d + \epsilon_{j+1})(\gamma + \epsilon_{j+1}^*)}}.\]
the exponent approaches $\frac{n \sigma}{d \gamma}$ as $j \to \infty$, as desired.
\paragraph*{Case 2: $I$ is essentially disjoint from the $T_i^{(j+1)}$.} In this case, we have the inequality
\[\mu_{j+1}(J) \leq 2 r^{n (L_{j+1} - L_j)} \mu_j(I)\]
but $N_{j+1}$ can be chosen sufficiently large so that $2 \mu_j(I) r^{-n L_j} < r^{-n \epsilon_{j+1} L_{j+1}}$, so we get the estimate
\[\mu_{j+1}(J) < r^{n L_{j+1}(1 - \epsilon_{j+1})}.\]
This estimate is at least as good as the desired estimate because $d \geq 1$ and $\sigma \leq \gamma$.

Now that we have proven the Frostman bound for the case where $I$ is a basic cube of the construction, it remains to show the Frostman estimate for arbitrary cubes $I$. As $\mu$ is a probability measure, we can restrict ourselves to the case for which $l(I) < r^{L_1}$. In particular, this means that there is some $j \geq 1$ such that $r^{L_{j+1}} \leq l(I) < r^{L_j}$. We will consider two cases: the case in which $r^{\lfloor (\sigma - 1000 \epsilon_{j+1}^*) N_{j+1}\rfloor} \leq l(I) < r^{L_j}$, and the complementary case in which $r^{L_{j+1}} \leq l(I) < r^{\lfloor (\sigma - 1000 \epsilon_{j+1}^*) N_{j+1} \rfloor}$.
\paragraph*{Case 1: $r^{\lfloor (\sigma - 1000 \epsilon_{j+1}^*) N_{j+1}\rfloor } \leq l(I) < r^{L_j}$.} Let $L$ be such that $l(I) = r^L$. In this case, $I \cap E$ can be covered by at most $Cr^{n(L - \lfloor (\sigma - 1000 \epsilon_{j+1}^*) N_{j+1} \rfloor)}$ cubes of sidelength $r^{\lfloor(\sigma - 1000 \epsilon_{j+1}^*) N_{j+1} \rfloor}$ occurring in step $j+1$ of the construction for some constant $C$ depending only on $n$ and $r$. Each of these cubes is known to have $\mu$-measure at most $\mu_j^* r^{n( \lfloor (\sigma - 1000 \epsilon_{j+1}^*) N_{j+1} \rfloor - L_j)}$, where $\mu_j^*$ is the maximum $\mu_j$-measure of any basic cube of sidelength $r^{L_j}$. Multiplying, we get that the $\mu$-measure of $I$ is at most $C \mu_j^* r^{n (L - L_j)}$.  But we have already established that $C \mu_j^* \leq r^{ \left(\frac{n \sigma}{d \gamma} - C' \epsilon_j \right)L_j}$ for some appropriate constant $C'$. Therefore, we get that  $\mu(I) \leq r^{n(L - L_j) + \left(\frac{n \sigma}{d \gamma} - C' \epsilon_j \right)L_j}$. We rearrange this expression to get $r^{\frac{n \sigma }{d \gamma } L + \left(n - \frac{n \sigma}{d \gamma} \right) (L- L_j) - C' \epsilon_j L_j}$. First, we observe that $\left(n - \frac{n \sigma}{d \gamma} \right) (L - L_j)$ is nonnegative, and thus the expression can only be made larger by removing this term. Second, because $L > L_j$, we have that $r^{-C' \epsilon_j L_j} \leq r^{-C' \epsilon_j L}$. Thus we get $\mu(I) \leq r^{\left(\frac{n \sigma}{d \gamma} - C' \epsilon_j \right) L}.$ The coefficient on $L$ approaches $\frac{n \sigma}{d \gamma}$ as $j \to \infty$, as desired.
\paragraph*{Case 2: $r^{L_{j+1}} \leq l(I) < r^{\lfloor (\sigma - 1000 \epsilon_{j+1}^*) N_{j+1} \rfloor}$.} By splitting $I$ into $O_n(1)$ parts, we may assume either $I \subset T_i^{(j+1)}$ for some $I$, or that $I$ is essentially disjoint from these sets. 

We begin with the case where $I$ is contained in some $T_{i}^{(j+1)}$. In this case, $I$ is contained in a union of at most $O_n(1)$ of the cubes of side length $r^{\lfloor(\sigma - 1000 \epsilon_{j+1}^*) N_{j+1} \rfloor}$ from Proposition \ref{NondegenerateProposition}. Therefore, up to an $O_n(1)$ loss, we can assume that $I$ is entirely contained in one of these cubes.

Let $L$ be such that $I$  has side length $r^L$. Then, by Part 2 of Proposition \ref{NondegenerateProposition}, $I$ intersects at most $\max(r^{n(L - \sigma N_{j+1})}, 1)$ of the cubes of side length $r^{L_{j+1}}$. 

If $L < \sigma N_{j+1}$, it follows that $\mu(I) \leq r^{(L - \sigma N_{j+1})n + L_{j+1} (\frac{n \sigma}{d \gamma} + C \epsilon_{j+1})}$. Using the relationship between $L_{j+1}$ and $N_{j+1}$, this is at most $r^{Ln -  C' \epsilon_{j+1} N_{j+1}}$ for some appropriate $C'$ depending on $n, d, \sigma$, and $\gamma$. But $L > (\sigma - 1000 \epsilon_{j+1}^*)N_{j+1}$, so this is no more than $r^{L(n - C'' \epsilon_{j+1})}$ for an appropriate $C''$ depending on $n, d, \sigma$, and $\gamma$. Notice that in this subcase we in fact get a bound that may be much smaller than $r^{L \frac{n \sigma}{d \gamma}}$.

If $L \geq \sigma N_{j+1}$, then $I$ intersects at most $1$ cube of side length $r^{L_{j+1}}$. This cube has measure at most $r^{(\frac{n \sigma}{d \gamma} - C \epsilon_{j+1} )L_{j+1}}$, and thus $\mu(I) \leq r^{(\frac{n \sigma}{d \gamma} - C \epsilon_{j+1})L}$ as desired since $L_{j+1} > L$.

We now consider the case where $I$ is essentially disjoint from the $T_i^{(j+1)}$. In this case, $I$ intersects at most $r^{n(L - L_{j+1})}$ of the cubes of sidelength $r^{L_{j+1}}$ that were retained as part of $E_{j+1}$. Each of these cubes has $\mu_{j+1}$-measure at most $r^{n L_{j+1} (1 - \epsilon_{j+1})}$. We multiply and conclude that 
\[\mu_{j+1}(I) \leq r^{nL - \epsilon_{j+1} L_{j+1}} \leq  r^{(n - C \epsilon)L}\]
because $L > (\sigma - 1000 \epsilon_{j+1}^*) N_{j+1}$, and thus $L > C^{-1}L_{j+1}$ for some appropriate constant $C$.
\section{Application to angle-avoiding sets}
M\'ath\'e \cite{Mathe17} established the following fact:
\begin{mythm}\label{MatheAngle}[Angle-Avoiding Sets, M\'ath\'e] Let $n \geq 2$, and let $\alpha \in (0, \pi)$ be such that $\cos^2(\alpha)$ is rational. There exists a compact set $E \subset \mathbb{R}^n$ of Hausdorff dimension $n/4$ such that $E$ does not contain three points forming an angle $\alpha$.
\end{mythm}
Theorem \ref{LandmarkTheorem} can be used to extend this result to all angles $\alpha$ for which $\cos \alpha$ is algebraic.
\begin{mythm}\label{AlgebraicAngle}[Angle-Avoiding Sets, Algebraic Case] Let $n \geq 2$, and let $\alpha \in (0, \pi)$ be such that $\cos \alpha$ is algebraic. Then there exists a compact set $E \subset \mathbb{R}^n$ of Hausdorff dimension $n/4$ such that $E$ does not contain three points forming an angle $\alpha$.
\end{mythm}
\begin{proof}
The proof is similar to the one in \cite{Mathe17}. We observe that three points $x, y, z \in \mathbb{R}^n$ form an angle $\alpha$ if they satisfy
\[(y - x) \cdot (z - x) = |y-x| |z-x| \cos \alpha.\]

We square both sides of this equation in order to turn the equation into a polynomial.

\[((y - x)\cdot (z - x))^2 = |y-x|^2 |z-x|^2 \cos^2 \alpha.\]

We then use the landmark system provided in Example \ref{AlgebraicExample} together with Theorem \ref{LandmarkTheorem} to conclude the desired result.
\end{proof}

M\'ath\'e \cite{Mathe17} proceeds to construct a set of Hausdorff dimension $n/8$ that avoids an arbitrary angle $\alpha$ by finding a polynomial with rational coefficients that vanishes on triples of points $(x,y,z)$ that form an angle $\alpha$. We now show that this example is typical.
\begin{mythm}
Let $p: \mathbb{R}^{nd} \to \mathbb{R} $ be a polynomial of degree $d$ whose coefficients lie in a $2$-dimensional vector space over $\mathbb{Q}$ of the form $\mathbb{Q} + \mathbb{Q} t$ for some number $t$. Then there exists a subset of $\mathbb{R}^n$ of Hausdorff dimension $\frac{n}{2d}$ that does not contain any $v$ distinct points $x_1, \ldots, x_v$ such that $p(x_1, \ldots, x_v) = 0$.
\end{mythm}
\begin{proof}
By multiplying by an appropriate integer, we can assume the coefficients of $p$ are in the finitely-generated free module $\mathbb{Z} + t \mathbb{Z}$.
 We assume the coefficients of $p$ are of the form $a + bt$, where $a$ and $b$ are integers. By Dirichlet's principle, there exist infinitely many pairs of integers $(r,q)$ such that $|t - r/q| \leq q^{-2}$. Therefore, the coefficients of $p$ are simultaneously approximable to degree $1$: $a + bt$ is within $b q^{-2}$ of the rational number $\frac{aq + br}{q}$. The same can also be said for all derivatives of the polynomial $p$. Using Theorem \ref{LandmarkTheorem} together with Example \ref{RationalApproximateExample} gives the desired result.
\end{proof}
This theorem can be extended in a trivial way:
\begin{mythm}
Let $p(x_1, \ldots, x_v)$ be a polynomial of degree $d$ whose coefficients lie in a $k$-dimensional vector space over $\mathbb{Q}$ of the form $\mathbb{Q} + \mathbb{Q} t_1 + \cdots + \mathbb{Q} t_{k-1} $ for some numbers $t_1, \ldots, t_{k-1}$. Then there exists a set of Hausdorff dimension $\frac{n}{dk}$ that does not contain any $v$ distinct points $x_1, \ldots, x_v$ such that $p(x_1, \ldots, x_v) = 0$.
\end{mythm}
\section{Polynomials in Non-archimedean Settings}
We can apply Theorem \ref{LandmarkTheorem} to conclude the following.

\begin{mycor}\label{LocalFieldPolynomial}
Let $p : K^{nv} \to K$ be a polynomial of degree $d$ on a non-archimedean local field with integer coefficients. If $K$ has characteristic $0$, or if $d < \text{char } K$, then there exists a subset of $K^n$ with Hausdorff dimension $n/d$ that does not contain any $v$ distinct points such that $p(x^{(1)}, \ldots, x^{(v)}) = 0$.
\end{mycor}
This follows from Theorem \ref{LandmarkTheorem} because any polynomial of degree $d$, where $d < \text{char } K$, will have a partial derivative of degree at most $d$ that is equal to a constant. The condition $d < \text{char } K$ is necessary for this observation to work: finite characteristic it is possible for a nonconstant polynomial of degree $\geq \text{char } K$ to have a derivative of zero. However, if $p$ is a degree $d$ polynomial where $d < \text{char } k$ than some appropriate $d$th partial derivative will be constant and nonzero and the assumptions of the theorem will therefore be met.

An important example of this occurs when $n = 1$ and $p(x, y, z) = x - 2y + z$. This is a polynomial that selects for three-term arithmetic progressions. In this case, Theorem \ref{LocalFieldPolynomial} states that there is a subset of $K$ with Hausdorff dimension $1$ that does not contain any nondegenerate $3$-term arithmetic progressions.

We focus especially on the case in which $K$ is the function field $\mathbb{F}_3((t))$. The unit ball, $\mathbb{F}_3[[t]]$, is isomorphic as a topological abelian group to the projective limit of the finite abelian groups $(\mathbb{Z}/3 \mathbb{Z})^n$, and thus the problem of finding large subsets of $\mathbb{F}_3[[t]]$ without $3$-term arithmetic progressions serves as a limiting case of the capset problem. The capset result states that for sufficiently large $n$, every subset of $(\mathbb{Z}/3 \mathbb{Z})^n$ with at least $2.756^n$ elements contains a $3$-term arithmetic progression \cite{EllenbergGijswijt17}. However, Corollary \ref{LocalFieldPolynomial} gives a subset of $\mathbb{F}_3[[t]]$ with Hausdorff dimension $1$ that does not contain any $3$-term arithmetic progressions, so a Hausdorff dimension analogue of the finite capset result does not hold in the limiting case. The author will consider the problem of a limiting capset result for \emph{Fourier dimension} in a future work.
\section{Acknowledgements}
This material is based upon work supported by the National Science Foundation under Award No. 1803086.
\bibliographystyle{myplain}
\bibliography{Master_Bibliography}
\end{document}